\documentstyle [twoside]{article}
\newtheorem{thm}{Theorem}[section]
\newtheorem{prop}[thm]{Proposition}
\newtheorem{cor}[thm]{Corollary}
\newtheorem{lem}[thm]{Lemma}
\newtheorem{remark}{Remark} 
\newtheorem{dfn}{Definition} 

\def\b0{{\bf 0}}
\def\ba{{\bf a}}
\def\bb{{\bf b}}
\def\be{{\bf e}}
\def\bg{{\bf g}}
\def\bh{{\bf h}}
\def\bi{{\bf i}}

\def\bz{{\bf z}}
\def\b1{{\bf 1}}

\def\bT{{\bf T}}

\def\bP{{\bf P}}
\def\bq{{\bf q}}

\def\bR{{\bf R}}

\def\bZ{{\bf Z}}
\def\bC{{\bf C}}
\def\bzeta{\zeta}

\def\Blat{\mbox{\it \raise2pt\hbox{"}\kern-2pt H}}
\def\lvup{\rlap{\ ${}^{\ell\atop{\hbox{${}^{\vee}$}}}$}\cdots}
\begin{document}
\begin{center}
{\center{\Large{\bf
Transposition mirror symmetry construction\\ 
and  period integrals
} }}

 \vspace{1pc}
{ \center{\large{ Susumu TANAB\'E }}}
\end{center}

\noindent
\begin{center}
 \begin{minipage}[t]{10.2cm}
{\sc Abstract.} {\em In this note we study several conditions
to be imposed on a mirror symmetry candidate 
to  the generic multi-quasihomogeneous Calabi-Yau
variety defined in the product of the quasihomogeneous projective spaces.
We propose several properties for a Calabi-Yau complete intersection variety
so that its period integrals can be expressed by means of quasihomogeneous 
weights of its mirror symmetry candidate as it has been suggested
by Berglund, Candelas et alii (Theorem ~\ref{thm31}).
As a corollary, we see certain duality between the monodromy data
and the Poincar\'e polynomials of the Euler characteristic for 
the pairs of our varieties (Theorem~\ref{thm432}).
}
 \end{minipage} \hfill
\end{center}
 \vspace{1pc}


{
\center{\section{ Introduction and Notations.}}
}

The transposition mirror construction has been proposed by P.Berglund and
T. H\"ubsch \cite{BH} as a trial to generalize so called
Greene-Plesser mirror construction
that comprises mirror pairs of Fermat type hypersurfaces. Later,
in the article \cite{Ber}, in relying on the transposition method,
the authors have proposed a natural hypothesis on the period integrals
associated to the complete intersection (CI) Calabi-Yau variety
$X$ that is supposed to be a mirror symmetry to the
generic multi-quasihomogeneous Calabi-Yau
variety $Y$ of codimension $\ell$ defined in the product of the
quasihomogeneous projective spaces
$\bP^{(\tau_1)}_{(g_1^{(1)},\ldots,
g_{\tau_1+1}^{(1)})} \times \ldots \times
\bP^{(\tau_k)}_{(g_1^{(k)},\ldots, g_{\tau_k+1}^{(k)})}$.
They mean under the notion of the mirror symmetry between $X$ and $Y$
an interchange between geometric symmetry $({\mathcal G}_X,$
${\mathcal G}_Y)$
  and quantum symmetry $({\mathcal Q}_Y,$ ${\mathcal Q}_X)$
groups of each varieties \cite{BK}
(See  Definition ~\ref{dfn11}, Theorem ~\ref{thm40} below).
That is to say, for the mirror symmetry pair $X$ and $Y$,
the following isomorphisms holds,
$${\mathcal Q}_X \cong {\mathcal G}_Y,$$
$${\mathcal Q}_Y \cong {\mathcal G}_X.$$
In this article we propose certain sufficient conditions
on $X$ and $Y$ so that
their hypothesis (\cite{Ber} \S 3.3) on the period integrals holds (Theorem ~\ref{thm31}).
Namely the  period integrals
defined on $X$ can be expressed by means of quasihomogeneous
weights of its mirror symmetry $Y$ and vice versa.
It will be shown that
our sufficient conditions entail the mirror symmetry between
$X$ and $Y$ in the above sense of \cite{BK} (see Theorem ~\ref{thm40}).


 \vspace{1.5pc}
\footnoterule

\footnotesize{AMS Subject Classification: 14J32 primary), 32S25, 32C65
(secondary).

Key words and phrases: mirror symmetry, hypergeometric functions,
monodromy.}

 \footnotesize{partially supported by Hokkaido University, ICTP (Trieste).}
\normalsize

\newpage
During the entire article we shall restrict ourselves to
the case $\ell =k.$

In accordance with the suggestion on the  mirror symmetry to the
generic multi-quasihomogeneous Calabi-Yau made in \cite{Ber},
we shall consider the following system of equations on
$\bT^n:=(\bC^\times)^n$ as defining equations of $X=X_\b1$,
$$X_s=\{x \in \bT^n; f_1(x)+ s_1(f_2(x)-1)= \cdots = f_{2k-1}(x)+
s_k (f_{2k}(x)-1)=0\},$$ with system $(1.1)$ below.
Here we use the notation
$s=(s_1, \cdots, s_k)\in \bT^k,$ $\b1=(1,1,\cdots,1) \in \bR^k.$
\begin{flushleft}
$$ f_1(x) = x^{\vec v^{(1)}_1} + \cdots + x^{\vec
v^{(1)}_{\tau_1}}, \leqno(1.1)$$
$$  f_2(x) =\prod_{j \in I^{(1)}} x_j  +1,$$
$$ \vdots,$$
$$ f_{2i-1}(x) = x^{\vec v^{(i)}_1} + \cdots + x^{\vec
v^{(i)}_{\tau_i}},$$
$$  f_{2i}(x) = \prod_{j \in I^{(i)}} x_j  +1,$$
$$ \vdots,$$
$$ f_{2k-1}(x) = x^{\vec v^{(k)}_1} + \cdots + x^{\vec
v^{(k)}_{\tau_k}},$$
$$  f_{2k}(x) = \prod_{j \in I^{(k)}} x_j+1.$$
\end{flushleft}
Here $I^{(j)}, 1 \leq j \leq k$ are
sets of indices which are complementary one another in such a way that
 $\cup_{q \in
[1,k]}I^{(q)}=\{1, \cdots, n \}$ and $I^{(q)}\cap I^{(q')}  =
\emptyset $ if $q \not = q'.$
From now on we shall make use  of the notations
$\tilde \tau_\nu := |I^{(\nu)}|$ and
$b^q :=\sum_{\nu=1}^q \tau_\nu.$
Additionally we suppose that
$$ \sum_{\nu=1}^k \tau_\nu =b^k = n.$$
The equation
$f_{2j-1}(x)$ (resp.$f_{2j}(x)$) is defined by the monomials
with powers
$\vec v^{(j)}_1,\cdots, \vec v^{(j)}_{\tau_j}
$ $\in \bZ^n$ (resp. $\vec v^{(j)}_{\tau_j+1} $ $\in \bZ^n$)
such that for the weight vector $\vec g^{(q)}=
(\overbrace{0,\cdots,0}^{b^{q-1}}, g_1^{(q)}, \cdots ,
g_{\tau_q}^{(q)}, \overbrace{0,\cdots,0}^{n-b^{q}})$ $\in
\bZ^{n}_{\geq 0}, $  $1 \leq q \leq k$
the following quasihomogeneiety condition holds,
$$ Q^{(q)}_j:=<\vec v^{(j)}_1, \vec g^{(q)}> = \cdots =
<\vec v^{(j)}_{\tau_j}, \vec g^{(q)}>=
<\vec v^{(j)}_{\tau_{j+1}},\vec g^{(q)}>
 \;\;\;, 1 \leq j \leq k. \leqno(1.2)$$
This means that the point
$\vec v^{(j)}_{\tau_j+1}$ belongs to the $(\tau_j-1)$-dimensional
hyperplane generated by
$\vec v^{(j)}_1, \vec
v^{(j)}_1,\cdots, \vec v^{(j)}_{\tau_j}.$
 The following condition shall be imposed if we
suppose that $X_s$ is a Calabi-Yau variety:
$$ \sum_{j=1}^kQ^{(q)}_j = \sum_{i=1}^{\tau_q}g_i^{(q)}
= <\vec g^{(q)},(\overbrace{0,\cdots,0}^{b^{q-1}},
\overbrace{1,\cdots,1}^{\tau_{q}},\overbrace{0,\cdots,0}^{n-b^{q}})>,
1 \leq q \leq k. \leqno(1.3)$$
We will impose this condition on $(1.1)$
in the further arguments.

In addition to that we assume that for each element
$ \lambda
\in Aut (X_j),$
of the group automorphism of the hypersurface
$X_j =\{x \in \bT^{n}; f_{2j-1}(x)=0 \}$ the relation
$(\lambda_{\ast}f_{2j-1})(x)= \lambda_{\ast}(x^{\vec
v^{(j)}_{\tau_j+1}})$ holds.
More precisely, every  $\lambda \in Aut (X_s)$
admits the following decomposition
$$ \lambda = \bq \cdot \bg \cdot \bh,$$
where $\bq \in {\mathcal Q}_X$, $\bg \in {\mathcal G}_X$ and
$\bh \in {\mathcal H}$ each of which is a non cyclic element of  $Aut (X_s)$.
The subgroup  ${\mathcal Q}_X \cong  \prod^{k}_{q=1}{\bf Z_{{\bar Q}^{(q)}}},$
with ${\bf {\bar Q}^{(q)}}= L.C.M.(Q^{(q)}_1, \cdots, Q^{(q)}_k) $ is a cyclic group generated by $k$ different cyclic actions
$\bq^{(\nu)}, 1 \leq \nu \leq k,$
corresponding to the quasihomogeneiety,
$$ \bq_\ast^{(\nu)}:(x_1, \cdots ,x_n) \rightarrow (x_1,
\cdots, x_{b^{\nu-1}}, e^{\frac{2 \pi i g_1^{(\nu)}}
{{\bf \bar Q}^{(\nu)}  }}x_{b^{\nu-1}+1}, \cdots,
e^{\frac{2 \pi i g_{\tau_\nu}^{(\nu)}}{{\bf \bar Q}^{(\nu)}}} x_{b^{\nu}}, x_{b^{\nu}+1}, \cdots, x_n).$$
The group (called geometric symmetry)
${\mathcal G}_X$ consists of elements $\bg \not \in {\mathcal Q}_X$
of the following form
$$ \bg_\ast:(x_1, \cdots ,x_n) \rightarrow (e^{\frac{2 \alpha_1 \pi i}{d}}x_1,
e^{\frac{2 \alpha_2 \pi i}{d}}x_2,\cdots, e^{\frac{2 \alpha_n \pi i}{d}}x_n),$$
for some $d>0$ and 
$(\alpha_1, \cdots, \alpha_n) \in \bZ^n.$ Here we remark that  
some of  $\alpha_i$ can be zero.
The group ${\mathcal H}$ is the non cyclic part of the group $Aut (X_s)$. 
\begin{dfn} In the following decomposition,
$$ Aut(X_s) \cong  {\mathcal Q}_X \times {\mathcal G}_X \times {\mathcal H},$$
we call ${\mathcal Q}_X$ (resp. ${\mathcal G}_X$) the quantum symmetry (resp.
geometric symmetry) of $X_s$.
\label{dfn11}
\end{dfn}
Further we apply so called
Cayley trick to $(1.1)$ to get a polynomial
$$ F(x,s,y) = \sum_{j=1}^{k} y_{2j-1}(f_{2j-1}(x)+s_j) +
\sum_{j=1}^{k} y_{2j}f_{2j}(x), \leqno(1.4)$$ with $L= n+3k$
terms. The procedures $(2.3)$, $(2.4)$ below explain
why we consider this polynomial $F(x,s,y)$
to calculate the period integrals associated to $X_s$.
One may consult \cite{Tan04} and \cite{Tan05} for more details
on the utility of the Cayley trick in the  calculus of
period integrals.

To manipulate the polynomial $F(x,s,y)$ we introduce the notation
 $a^\nu:=\tau_1 + \cdots +\tau_{\nu}+3\nu = b^\nu+ 3\nu.$
In particular, the $a^{i-1}+1-$th term of $F(x,s,y)$
corresponds to $y_{2i-1}x^{\vec v^{(i)}_1}$ and the $a^{i}-3=
(a^{i-1}+\tau_i)-$th term  - $y_{2i-1}x^{\vec
v^{(i)}_{\tau_i}}.$ The $(a^{i}-1)-$th term  - $y_{2i} \prod_{j
\in I^{(i)}} x_j.$ The$(a^{i}-2)-$th term - $y_{2i}s_i.$ The
$a^{i}-$th term   - $y_{2i}$.

From the polynomial $ F(x,s,y)$ we construct a matrix $\sf L$
that consists of the row vectors $r-$th term of which corresponds
to the power of the $r-$th monomial term present in $ F(x,s,y).$
For example, the row vector $\bar v_{q}^{(\nu)}$
corresponds to the row number $a^{q-1}+ \nu$ of  the matrix $\sf L$,
$$ \bar
v_{q}^{(\nu)}= (\vec v_{q}^{(\nu)},\overbrace {0,
\cdots,0}^{2\nu-2},1, \overbrace{0, \cdots, 0}^{3k-2\nu
+1}),$$ $1 \leq \nu \leq k, 1 \leq q \leq \tau_\nu.$

Next we look at the system of linear equations
$\Xi = (\xi_1, \cdots, \xi_L),$
$$ ^t{\sf L}\cdot {\Xi} = ^t(1,\cdots,1, z_1, \cdots
,z_k),\leqno(1.5)$$ that is equivalent to the relation,
$${\Xi}=^t{\sf L}^{-1} \cdot ^t(1,\cdots,1, z_1, \cdots ,z_k).
\leqno(1.6)$$ Thus we get linear functions
$(\xi_1(z), \cdots, \xi_L(z))$ that will be later denoted by
$\left({\mathcal L}_1(0,z,0), \cdots, {\mathcal L}_L(0,z,0) \right)$
in the notation $(2.5)$ below.
These linear functions will play essential role
in the calculus of the period integrals.

We define a $(n \times k)-$ matrix $V^\Lambda$ as follows:
$$ V^\Lambda := ( ^t{\vec v}_{\tau_1+1}^{(1)}, \cdots,
^t{\vec v}_{\tau_k+1}^{(k)}), \leqno(1.7)$$ where
${\vec v}_{\tau_q+1}^{(q)}$ is a  $n-$row vector that corresponds to
$supp(f_{2q}) \setminus \{0\}.$
By virtue of the
quasihomogeneiety $(1.2)$ we can define a  $k\times
k$ matrix as follows:
$$  \hat{Q} :=
\left (\begin {array}{ccc}
Q^{(1)}_1&\cdots&Q^{(k)}_1\\
\vdots&\vdots&\vdots\\
Q^{(1)}_k&\cdots&Q^{(k)}_k\\
\end {array}\right )
=  ^tV^\Lambda \cdot \left(
^t\vec g^{(1)},\cdots, ^t\vec g^{(k)}\right ).
\leqno(1.8)$$
For the simplicity of the formulation, we
will use a diagonal matrix
$$ G=diag(g^{(1)}_1,\cdots, g^{(1)}_{\tau_1}, g^{(2)}_1,
\cdots, g^{(2)}_{\tau_2}, \cdots, g^{(k)}_1,\cdots,
g^{(k)}_{\tau_k} ). \leqno(1.9)$$
We introduce a $n\times n$ matrix:
$$  {\sf L}_\Lambda :=\left (\begin {array}{c}
\vec v_1^{(1)}\\\vec v_2^{(1)}\\
\vdots\\
\vec v_{\tau_1}^{(1)}\\
\vdots\\
 \vec v_{\tau_k}^{(k)}\\
\end {array}\right ). \leqno(1.10)$$

We construct a matrix $^T\sf L$ constructed from the
transposed matrix $^t\sf L$ after some proper permutations of the rows
and columns
such that each row of
$^T \sf L$ corresponds to a vertex of a polynomial
$$ ^TF(x,s,y)=\sum_{j=1}^{k} y_{2j-1}(^Tf_{2j-1}(x)+s_j) +
\sum_{j=1}^{k} y_{2j}\;^Tf_{2j}(x), \leqno(1.4)^T$$
for the polynomials,
 \begin{flushleft}
$$ ^Tf_{2q-1}(x) = x^{^T\vec v^{(q)}_1} + \cdots + x^{^T\vec
   v^{(q)}_{\tilde \tau_q}}, \leqno (1.1)^T$$
$$ ^Tf_{2q}(x) =\prod_{\ell \in ^TI^{(q)}} x_\ell  +1,\;\;,
\; 1 \leq q \leq k.  $$
\end{flushleft}
Here we impose the condition
$$ \{\tau_1, \cdots, \tau_k\}=
\{ |\;^TI^{(1)}|, \cdots,  |\;^TI^{(k)}| \}
= \{\tilde \tau_1, \cdots, \tilde \tau_k\}
= \{ |I^{(1)}|, \cdots,  |I^{(k)}| \}. \leqno(1.11)$$

Thus we can define an one to one mapping
$$ \nu:[1,k] \rightarrow [1,k],$$
such that $ |^TI^{(\nu(j))}|=
\tau_{\nu(j)} =\tilde \tau_j = |I^{(j)}|$.
Further we impose a condition
$\nu^2=id$ so that $\;^T(^TF(x,s,y))= F(x,s,y)  $ holds.
In a way parallel to $(1.7) - (1.10)$,
we can define the weight
system
$$ (^T\vec g^{(1)}, \cdots , ^T\vec g^{(k)}), $$
$$^T{Q}_j^{(q)}= \langle ^T\vec v^{(j)}_r,  ^T\vec g^{(q)}\rangle,\; 1
\leq r \leq \tilde \tau_q, 1 \leq j,q \leq k. \leqno(1.2)^T$$
We impose a condition necessary for Calabi-Yau property of $Y$
$$ \sum_{j=1}^k\;^TQ^{(q)}_j = \sum_{i=1}^{\tilde \tau_q}\;^Tg_i^{(q)},
1 \leq q \leq k. \leqno(1.3)^T$$
It is easy to see that the equations of  $\;^T(1.1)$
define a  CI in
$\bP^{(\tilde \tau_1)}_{(^Tg_1^{(1)},\ldots,
^Tg_{\tilde \tau_1+1}^{(1)})} \times
\ldots \times \bP^{(\tilde \tau_k)}_{(^Tg_1^{(k)},\ldots,
^Tg_{\tilde \tau_k+1}^{(k)})}$
with $ ^Tg_{\tilde \tau_q+1}^{(q)}= ^T{Q}_q^{(q)}.$

In analogy with the Definition ~\ref{dfn11} we define the following decomposition
$$ Aut(Y_s) \cong  {\mathcal Q}_Y \times {\mathcal G}_Y \times \;^T{\mathcal H},$$
for
$$Y_s=\{x \in \bT^n; \;^Tf_1(x)+ s_1(\;^Tf_2(x)-1)= \cdots = \;^Tf_{2k-1}(x)+
s_k (\;^Tf_{2k}(x)-1)=0\}.$$
The quantum symmetry ${\mathcal Q}_Y \cong  
\prod^{k}_{q=1}{\bf Z_{\;^T{\bar Q}^{(q)}}},$
with ${\bf\;^T{\bar Q}^{(q)}}= L.C.M.(\;^TQ^{(q)}_1, \cdots, \;^TQ^{(q)}_k) $ 
is a cyclic group corresponding to the quasihomogeneiety  and  
${\mathcal G}_Y $ the remaining cyclic part called geometric symmetry.
The group $\;^T{\mathcal H}$ is the remaining non cyclic part of $ Aut(Y_s).$

We introduce matrices analogous to the case of
$(1.1),$ $$^TV^\Lambda := (
^t(^T\vec v_{\tau_1+1}^{(\nu(1))}), ^t(^T \vec v_{\tau_2+1}^{(\nu(2))}),
\cdots, ^t(^T \vec v_{\tau_k+1}^{(\nu(k))})), \leqno(1.7)^T$$ where
$^T\vec v_{\tau_q+1}^{(\nu(q))}$ is a  $n-$row
vector which corresponds to
$supp(^Tf_{2\nu(q)}) \setminus \{0\}.$
More precisely, the $j-$th column of the matrix  $^TV^\Lambda$
equals to
$$ ^t(^T\vec v_{\tau_j+1}^{(\nu(j))})=
^t(\overbrace{0, \cdots,0}^{b^{\nu(j)-1}}, 
\overbrace{1, \cdots,1}^{\tau_{\nu(j)}},
\overbrace{0, \cdots,0}^{n-b^{\nu(j)}}), \tau_{\nu(j)}=\tilde {\tau_j},$$
here $b^{\nu(j)} := \sum^{\nu(j)}_{r=1} \tau_r.$
$$  \;^T{\hat Q} :=
\left (\begin {array}{ccc}
^TQ^{(1)}_{\nu(1)}&\cdots&^TQ^{(k)}_{\nu(1)}\\
\vdots&\vdots&\vdots\\
^TQ^{(1)}_{\nu(k)}&\cdots&^TQ^{(k)}_{\nu(k)}\\
\end {array}\right )
= ^t(^TV^\Lambda) \cdot
\left (^t(^T\vec g^{(1)}),\cdots,^t(^T \vec g^{(k)})\right ). \leqno (1.8)^T$$
$$ ^TG=diag
(^Tg^{(1)}_1,\cdots, ^Tg^{(1)}_{\tilde \tau_1}, ^Tg^{(2)}_1, \cdots,
^Tg^{(2)}_{\tilde \tau_2}, \cdots, ^Tg^{(k)}_1,\cdots,
\;^Tg^{(k)}_{\tilde \tau_k}).    \leqno(1.9)^T$$
$$  ^T{\sf L}_\Lambda :=\left (\begin {array}{c}
^T\vec v_1^{(\nu(1))}\\
^T\vec v_2^{(\nu(1))}\\
\vdots\\
^T\vec v_{\tilde \tau_{\nu(1)}}^{(\nu(1))}\\
\vdots\\
^T\vec v_{\tilde \tau_{\nu(k)}}^{(\nu(k))}\\
\end {array}\right )
=\left (\begin {array}{c}
^T\vec v_1^{(\nu(1))}\\
^T\vec v_2^{(\nu(1))}\\
\vdots\\
^T\vec v_{\tau_1}^{(\nu(1))}\\
\vdots\\
^T\vec v_{\tau_{k}}^{(\nu(k))}\\
\end {array}\right ).  \leqno(1.10)^T $$
In an analogous way to $(1.6)$,
we introduce the linear functions
$^T\Xi := (^T\xi_1(z), \cdots,^T\xi_L(z)),$
defined by the relation,
$$^T{\Xi}
= \;^t({^T\sf L})^{-1} \cdot ^t(1,\cdots,1, z_1, \cdots ,z_k).\leqno(1.6)^T$$
Finally we remark that due to the property
$\nu^2=id$ and the definition $(1.7)^T$, there exists a permutation matrix
$\lambda \in SL(n, \bZ)$ such that
$$ \lambda \cdot V^\Lambda = \;^TV^\Lambda, \lambda \cdot {^t\sf L_\Lambda}= {^T\sf L_\Lambda}.
\leqno(1.12)$$

{
\center{\section{ Mellin transform of period integrals }}
}

In this section we review the results on the period integrals 
to be used for the verification of the hypothesis in \cite{Ber}
in the subsequent section. See for the detail of proofs
\cite{Tan04}, \cite{Tan05}.

Let us consider the Leray's coboundary (see \cite{Vas}) to define the
period integral that is equivalent to the period integral of the variety 
$X_s$,
$\gamma \subset H_{n}({\bf T}^{n} \setminus
\cup_{\ell=1}^k \{x \in \bT^n: f_{2\ell-1}(X)+s_\ell=0\}
\cup_{\ell=1}^k \{x \in \bT^n: f_{2\ell}(X)=0\})$ such that
$\Re (f_{2\ell-1}(X) + s_\ell)|_{\gamma} <0$,
$\Re f_{2\ell}(X)|_{\gamma} <0$.
Further on central object of our study is the following fibre integral,
$$  I_{x^{\bi}, \gamma}^{\zeta}(s)=
\int_{\gamma}
 (f_1(x) +s_1)^{-\zeta_1-1}(f_2(x))^{-\zeta_2-1}
\cdots (f_{2k-1}(x) +s_k)^{-\zeta_{2k-1}-1}
f_{2k}(x)^{-\zeta_{2k}-1} x^{\bi+\b1} \frac{dx}{x^{\b1}},
\leqno(2.1)$$ and its Mellin transform,
$$ M_{{\bi},\gamma}^\zeta ({\bz}):=\int_\Pi s^{\bz }
I_{x^{\bi}, \gamma}^{\zeta} (s)
\frac{ds}{s^{\b1}},
\leqno(2.2)$$ 
for  certain cycle $\Pi$ homologous to ${\bf R}^{k}$
which avoids the singular loci of 
$I_{x^{\bi}, \gamma}^{\bzeta} (s)$  (cf. \cite{Tsikh}).
Thus the fibre integral 
$I_{x^{\bi}, \gamma}^{\zeta}(s)$ is a ramified function on the torus
$\bT^k.$
We introduce the  notation $\gamma^\Pi:= \cup_{(s)\in
\Pi}((s), \gamma).$ One shall not confuse it with the 
thimble of Lefschetz, because  $\gamma^\Pi$ is rather a tube without thimble.
It is useful to understand the calculus of the Mellin
transform
in connection with the notion of the generalized
HGF in the sense of Mellin-Barnes-Pincherle \cite{Nor}.
After this formulation, the classical HGF of Gauss can be expressed 
by means of the integral,
$$ _2 F_1(\alpha,\beta,\gamma|s)=
\frac {1}{2\pi i } \int_{z_0 - i\infty}^{  z_0 + i\infty }(-s)^z
\frac {\Gamma ( z+ \alpha)\Gamma ( z+ \beta )\Gamma ( -z)}{\Gamma
( z+ \gamma)} dz , \;\; - \Re  \alpha, - \Re  \beta < z_0.$$
We can introduce $(w_1'', \cdots, w_{2k}'')$ natural quasihomogeneous weight
of $(y_1, \cdots y_k)$ so that $F(x,0,y)$ of $(1.4)$ gets the quasihomogeneous
zero weight with respet to the variables $(x,y).$
Next we modify the Mellin transform
$$
M_{{\bi},\gamma}^\zeta ({\bz})
= c(\zeta) \int_{S^{k-1}_+(w'') \times \gamma^\Pi} \frac{x^{\bi}
\omega^\zeta s^{\bz -\b1} dx \wedge \Omega_0(\omega) 
\wedge 
ds}
{(\omega_1(f_1(x) +s_1)+\cdots +\omega_{2k}(f_{2k}(x)))^ {\zeta_1+
\cdots +\zeta_{2k} +2k}} $$ $$= c(\zeta)\int_{\bR_+} \sigma^{\zeta_1+ \cdots
+\zeta_{2k} +2k} \frac{d \sigma}{\sigma} \int_{S^{2k-1}_+(w'')}
\omega^\zeta \Omega_0(\omega) \int_{\gamma}x^{\bi} dx \int_{\Pi}
s^{\bz} e^{\sigma(\omega_1(f_1(X) +s_1)+\cdots+ 
\omega_{2k}(f_{2k}(x))} \frac{ds}{s^\b1},$$
with $ c(\zeta) = \frac{\Gamma(\zeta_1+\cdots +\zeta_{2k}+{2k})}
{\Gamma(\zeta_1+1)\cdots \Gamma(\zeta_{2k}+1)}.$ Here we made 
use of notations
$S^{{2k}-1}_+(w'') =\{(\omega_1, \cdots,
\omega_{2k}): \omega_1^{|\frac{\bf w''}{w_1''}|}+\cdots
+\omega_{2k}^{|\frac{\bf w''}{w_{2k}''}|}=1, \omega_\ell >0\;\;$ for 
all  $\ell,$ ${\bf w''}$ $=$ $\prod_{1 \leq i \leq {2k}}w''_i \}$ and
$\Omega_0(\omega)$ the $({2k}-1)$ volume form on $S^{{2k}-1}_+
(w''),$
$$\Omega_0(\omega) = \sum_{\ell=1}^{2k} (-1)^\ell w''_\ell \omega_\ell d\omega_1\wedge \lvup \wedge
d\omega_{2k}.$$
In the above transformation we used a classical interpretation of
Dirac's delta function as a residue:
$$\int_{\gamma} \int_{\bR_+}
e^{ y_j(f_j(x) + s_j)} y_j^{\zeta_j} dy_j \wedge dx
= \Gamma(\zeta_j+1) \int_{\gamma}
 (f_j(x) + s_j)^{-\zeta_j-1} dx. $$

We will rewrite, up to constant multiplication, the expression
obtained as a modification of $M_{{\bi},\gamma}^\zeta ({\bz})$
into the following form,
$$\int_{(\bR_+)^{{2k}} \times \gamma^\Pi}
e^{\Psi(T)} x^{\bi+\b1} y^{\zeta+\b1}s^{\bz} \frac{dx}{x^\b1}
\wedge \frac{dy}{y^\b1} \wedge \frac{ds}{s^\b1}    
$$
where
$$\Psi(T) = T_1(X,s,y) + \cdots + T_L(x,s,y)= F(x,s,y),  \leqno(2.3)$$
in which each term $T_i(x,s,y)$ stands for a monomial in 
variables
$(x,s,y)$ of the phase function $(1.4).$
We transform the above integral into the following form,
$$
\int_{({\bR_+})^{2k} \times  \gamma^\Pi } e^{\Psi(T(x,s,y))}
x^{\bi+\b1 } s^\bz y^{{\bzeta}+\b1} \frac{dx}{x^\b1} \wedge
\frac{dy}{y^\b1} \wedge \frac{ds}{s^\b1} \leqno(2.4)$$
$$
= (det {\sf L})^{-1}\int_{{\sl L}_\ast ({ \bR_+}^{2k}
\times \gamma^\Pi  )}
e^{\sum_{a \in I}T_a} \prod_{a\in I}
T_a^{{\mathcal L}_a({\bi, \bz, \bzeta})}
\bigwedge_{a \in I} \frac{dT_a}{T_a}  $$
$$
= (-1)^{\zeta_1+\cdots +\zeta_{2k} +{2k}}
(det {\sf L})^{-1}\int_{-{\sl L}_\ast ({ \bR_+}^{2k}
\times \gamma^\Pi  )}
e^{-\sum_{a \in I}T_a} \prod_{a\in I}
T_a^{{\mathcal L}_a({\bi, \bz, \bzeta})}
\bigwedge_{a \in I} \frac{dT_a}{T_a}.  $$
Here ${\sl L}_\ast ({\bR_+}^{2k}
\times \gamma^\Pi)$ means the image of the chain in
$\bC^{M}_X\times \bC^{k}_s\times \bC^{{2k}}_y$
into that in  $\bC^{L}_T$
induced by the  transformation $(2.3).$
We define $-{\sl L}_\ast ({ \bR_+}^{2k}
\times \gamma^\Pi  )$ $=\{(-T_1, \cdots, -T_L)$ $\in \bC^L;$ 
$(T_1, \cdots,T_L)$ $\in
{\sl L}_\ast ({ \bR_+}^{2k}
\times \gamma^\Pi  ), \Re T_a >0, a \in [1,L] \}.$
The second equality of $(2.4)$ follows from Proposition 2.1, 3) below
that can be proven in a way independent of the argument to derive $(2.4).$
We will denote the set of columns and rows of the matrix 
$\sf L$ by  $I,$
$$ I:= \{1, \cdots, L\}.$$ Here we remember the  condition
$L=n+3k$ imposed on $(1.4)$.

The following notion helps us to formulate the result in a compact manner.
\begin{dfn}
A meromorphic function $g(\bz)$
is called $\Delta-$periodic for  
$\Delta \in \bZ_{>0},$ if
$$g(\bz)= h(e^{2 \pi \sqrt -1
\frac{z_1}{\Delta}}, \cdots,
e^{2 \pi \sqrt -1
\frac{z_k}{\Delta}} ),$$
for some rational function $h(\zeta_1, \cdots, \zeta_{k}).$
\label{dfn21}
\end{dfn}

For the CI  $(1.1)$ (i.e. we can construct $F(x,s,y)$ for which
the matrix $\sf L$ is non-degenerate), we have the following statement.
\begin{prop}
1)For any cycle $\Pi \in H_{k}({\bf T}^{k} \setminus S.S.
I_{x^{\bi}, \gamma}^{\bzeta} (s))$ 
the Mellin transform $(2.1)$
can be represented as a product of $\Gamma-$ function factors
up to a $\Delta-$periodic function factor $g(\bz)$,
$$  M_{{\bi}, \gamma}^\zeta (\bz)= g(\bz) \prod_{a \in I}
\Gamma\bigl( {\mathcal L}_a({\bi, \bz, \bzeta})\bigr),$$ with
$${\mathcal L}_a({\bi, \bz, \bzeta} ) =
\frac{\left(\sum_{j=1}^n A_j^a (i_j+1)
+\sum_{\ell=1}^k  B_\ell^a z_\ell +\sum_{\ell=1}^{2k}
D_\ell^a(\zeta_\ell+1)\right)}{\Delta}, a \in I. \leqno(2.5)$$
Here the following matrix $\Delta^{-1}{\sf T} = (\sf L)^{-1}$ 
has integer elements,
$$^t{\sf T}=(A_1^a,\cdots,
B_1^a, \cdots, B_k^a ,D_1^a
, \cdots, D_{2k}^a )_{1 \leq
a \leq L}, \leqno(2.6)$$
with $G.C.D.(A_1^a,\cdots,
A_{n}^a, B_1^a, \cdots, B_k^a ,D_1^a
, \cdots, D_{2k}^a )=1,$
for all $ a \in I.$
In this way $\Delta >0$ is uniquely determined.

The coefficients of  (2.5) satisfy the following properties for each index
 $a \in I$ :

$\bf a$
Either
${\mathcal L}_a({\bi, \bz, \bzeta} ) = \frac{\Delta}{\Delta}z_\ell,$
i.e. $A_1^a=\cdots=
A_{n}^a=0,$ $B_1^a = \cdots\lvup = B_k^a=0,$ $B_\ell^a=1.$

$\bf b$
Or
${\mathcal L}_a({\bi, \bz, \bzeta} ) = \frac{\Delta}{\Delta}(\zeta_{2\ell-1}
+\zeta_{2\ell}-z_\ell),$

$\bf c$
Or
$${\mathcal L}_a({\bi, \bz, \bzeta} )=
\frac{\sum_{j=1}^n A_j^a (i_j+1)
+\sum_{\ell=1}^k B_\ell^a (z_\ell -\zeta_{2\ell-1}-1) }{\Delta}$$
2) For each fixed index $1 \leq \ell \leq n, 1 \leq q \leq k,$ 
the following equalities take place:
$$\sum_{a \in I} A_\ell^a =0,\; \sum_{a \in I}
B_q^a =0. \leqno(2.7)$$

3) The following relation holds among the linear functions
 ${\mathcal L}_a,$ $a \in I$:
$$ \sum_{a \in I }{\mathcal L}_a(\bi, \bz, \zeta)= \zeta_1 +\cdots +
\zeta_{2k} +2k.$$
\label{prop21}

\end{prop}

In the sequel, especially from $(3.2)$ of the next section, 
we will make use of the notation as follows,
$$ (w_1^{(a)}, w_2^{(a)},\cdots,w_{n}^{(a)},p_1^{(a)},\cdots,p_{k}^{(a)},
q_1^{(a)},\cdots,q_{2k}^{(a)})=
(\frac{A_1^a}{\Delta},\cdots,
\frac{A_n^a}{\Delta},  \frac{B_1^a}{\Delta}
, \cdots,  \frac{B_k^a}{\Delta},  \frac{D_1^a}{\Delta}
, \cdots, \frac{D_{2k}^a}{\Delta})_{,1 \leq
a \leq L}. \leqno(2.8)$$

\par 
In view of the Proposition ~\ref{prop21},  we introduce the subsets
of indices $a \in I= \{1,2, \cdots, L\}$ as follows.
\begin{dfn}
The subset $I^+_q \subset I$ (resp. $I^-_q, I^0_q$)
consists of the indices
$a$ such that the coefficient $B^a_q$ of
${\mathcal L}_a(\bi, \bz, \bzeta)$ (2.5)
is positive (resp. negative, zero).  
\label{dfn2}
\end{dfn}

From this proposition we get the following.

\begin{cor}
The integral $I_{x^{\bi}, \gamma}^{\bzeta} (s)$ satisfies the hypergeometric
system of Horn type as follows:
$$ L_{q}(\vartheta_s, s)
I_{x^{\bi}, \gamma}^{\bzeta} (s):=
\Bigl[P_{q, \bi}( \vartheta_s, \zeta)
-s_q^\Delta Q_{q, \bi}( \vartheta_s, \zeta)\Bigr]
I_{x^{\bi}, \gamma}^{\bzeta} (s)=0,1 \leq q \leq k \leqno(2.9)_1$$
with
$$P_{q, \bi}(\vartheta_s, \zeta)=
 \prod_{a\in I^+_q} \prod_{j=0}^{B_q^a-1}
\bigl({\mathcal L}_a(\bi,-\vartheta_{s},
\zeta)+j\bigr), \leqno(2.9)_2$$ 
$$Q_{q, \bi}( \vartheta_s, \zeta) = \prod_{\bar{a}\in
I^-_q}\prod_{j=0}^{-B_q^{\bar a}-1} \bigl({\mathcal
L}_{\bar{a}}(\bi,-\vartheta_{s}, \zeta)+j \bigr)
,\leqno(2.9)_3$$
where $I^+_q, I^-_q,  1 \leq q \leq k $
are the sets of indices defined in  Definition ~\ref{dfn2}.

The degree of two operators  $P_{q,
\bi}(\vartheta_s, \zeta),$ 
$Q_{q, \bi}( \vartheta_s, \zeta)$ are equal.
Namely,
$$ deg\;P_{q,
\bi}(\vartheta_s, \zeta) =  \sum_{a \in I^+_q} B_q^a 
= - \sum_{\bar a \in I^-_q} B_q^{\bar a} =deg\;Q_{q,
\bi}(\vartheta_s, \zeta) .
\leqno(2.10)$$ 

\label{cor2}
\end{cor}

{
\center{\section{ Hypothesis by Berglund, Candelas et alii}}
}

In this section, we apply the results from \S 2 on the period integrals
to a class of Calabi-Yau varieties $(1.1)$ studied in the framework
of Landau-Ginzburg vacua theory.

Our aim is to find out sufficient conditions so that the (Mellin transform of)
the period integrals on $X_s$ can be expressed by means of quasihomogeneous
weight data of $Y.$ The main theorem of this section is
Theorem~\ref{thm31}.

Before proceeding to the proof of the Theorem~\ref{thm31},
we write down concretely the matrix $\sf L$ 
in taking  $(1.1), (1.4)$ into account, 
and we have the following 
row vectors of $\sf L$:
$$ \bar
v_{q}^{(\nu)}= (v_{q,1}^{(\nu)}, v_{q,2}^{(\nu)},\cdots,
v_{q,\tau_1}^{(\nu)}, \cdots, v_{q,n}^{(\nu)},\overbrace {0,
\cdots,0}^{2\nu-2},1, \overbrace{0, \cdots, 0}^{3k-2\nu
+1}),\leqno(3.1)$$ $1 \leq \nu \leq k, 1 \leq q \leq \tau_\nu.$
In using the notations $(2.8), $ $\S 2$ 
for the column vectors of the matrix
${\sf L}^{-1},$ one can deduce the following system for each 
$\nu$ and $q$.
$$ \leqno(3.2)\begin{array}{ccc} v_{q,1}^{(\nu)} w_1^{(a^\nu)} +
v_{q,2}^{(\nu)}w_2^{(a^\nu)} + \cdots+ v_{q,n}^{(\nu)}
w_{n}^{(a^\nu)}&=&-1\;\; for\;\; a^{\nu-1}+1 \leq q \leq
a^{\nu-1}+\tau_\nu,
\\ &=&0 \;\;\rm{otherwise}.
\end{array} $$
$$
v_{1,j}^{(1)} p_q^{(1)} + v_{2,j}^{(1)}p_q^{(2)} + \cdots+
v_{\tau_1,j}^{(1)} p_q^{(\tau_1)} + v_{1,j}^{(2)}p_q^{(\tau_1+4)}
+ \cdots+ v_{\tau_2,j}^{(2)} p_{q}^{(\tau_1+\tau_2+3)}+ \cdots +
v_{\tau_k,j}^{(k)} p_q^{(a^k-3)} \leqno(3.3) $$
$$\begin{array}{ccc}
&=&-1\;\; for
\;\;j \in I^{(q)},\\ 
&=&0  \;\; for \;\; j \in I^{(r)}, r \not = q.
\end{array} $$
Here we shall remark that the system
$(3.2)$ for a fixed  $\nu$
(resp. $(3.3)$ for a fixed  $q$) consists of $n-$ linear 
independent equations with respect to unknowns
$ \{w_1^{(a^{\nu})}, \cdots, w_n^{(a^{\nu})} \},$ (resp. $ \{
p_q^{(j)};j \in I_\Lambda \}$ ).
Here we made use of the notation
$I_\Lambda$ the set of indices 
$ \{1, \cdots, L\} \setminus
\cup_{\nu=1}^k\{a^{\nu}-2,a^{\nu}-1,a^{\nu} \}.$ 
One sees other necessary conditions on
$w_j^{(a^{\nu-1})}$, $$
\sum_{j \in I^{(\nu)}} w_j^{(a^{\nu}-1)}=1 , \leqno(3.4) $$
$$ \sum_{j \in
I^{(\nu)}}w_j^{(a^{\nu})}=-1 . \leqno(3.5) $$ 
This can be seen from the fact that the product of the
$(a^{\nu}-1)-$th column of the matrix 
${\sf L}^{-1}$ with the
$(a^{\nu}-1)-$th row of  ${\sf L}$ is equal to $1.$ 
On the other hand
$n+2\nu-$th column of ${\sf L}$ with the $j-$th row
($1 \leq j \leq n+2$) of  ${\sf L}^{-1}$ is equal to  $0$
which entails $w_j^{(a^{\nu}-1)} +w_j^{(a^{\nu})} =0.$
In addition to that, we see 
$$
\sum_{j \in I^{(\nu)}} w_j^{(a^{\nu'}-1)}= \sum_{j \in
I^{(\nu)}}w_j^{(a^{\nu'})}=0,$$
for $\nu' \not = \nu.$

In this situation we deduce the following system from
$(3.3)$, 
$$
\left (\begin {array}{c}
\vec P_1\\
\vec P_2\\
\vdots\\
\vec P_k\\
\end {array}\right )
\cdot {\sf L}_\Lambda =- \;^tV^\Lambda, \leqno(3.6)$$
with $\vec P_q =(p_q^{(1)},p_q^{(2)}, \cdots, p_q^{(\tau_1)}, \cdots,
p_q^{(a^k-3)} ).$
If we denote by  $\check{\sf L}_j$ the $j-$th column vector of the matrix
${\sf L}_\Lambda,$ we have the following equation 
derived directly from $(3.6)$:
$$ \langle  z_1 \vec P_1 + z_2 \vec P_1 + \cdots + z_k
\vec P_k, \check{\sf L}_j
\rangle+z_q=0\;\; {\rm{if}} \;\; j \in I^{(q)}. \leqno(3.7)$$
 In view of $(1.12)$, we get
$$ ({^T\sf L}_\Lambda)^{-1}\cdot ^TV^\Lambda
= ({^t\sf L}_\Lambda)^{-1}\cdot V^\Lambda, \leqno(3.8)$$
that yields
$$
-\left(\begin {array}{c}
^T\vec g^{(1)}\\
\vdots\\
^T\vec g^{(k)}\\
\end {array} \right )
\cdot ^TG^{-1}\cdot ({^T\sf L}_\Lambda)^{-1}
\cdot ^TV^\Lambda \leqno(3.9)$$
$$=\left (\begin {array}{cccc} 
p_1^{(1)}+ \cdots +p_1^{(\tilde \tau_1)} &p_2^{(1)}+ \cdots +p_2^{(\tilde \tau_1)}& \cdots & 
p_k^{(1)}+ \cdots +p_k^{(\tilde \tau_1)}\\
p_1^{(\tilde a_1+1)}+ \cdots +p_1^{(\tilde a_1+\tilde \tau_2)} &p_2^{(\tilde a_1+1)}+ \cdots 
+p_2^{(\tilde a_1+\tilde \tau_2)}& \cdots & 
p_k^{(\tilde a_1+1)}+ \cdots +p_k^{(\tilde a_1+\tilde \tau_2)}\\
\vdots&\vdots&\vdots &\vdots \\
p_1^{(\tilde a_{k-1}+1)}+ \cdots +p_1^{(n+3k-3)} &
p_2^{(\tilde a_{k-1}+1)}+ \cdots 
+p_2^{(n+3k-3)}& \cdots & 
p_k^{(\tilde a_{k-1}+1)}+ \cdots +p_k^{(n+3k-3)}\\
\end {array}\right ),$$
where we used the notation $\tilde a_q := \sum_{\nu=1}^q(\tilde \tau_\nu+3).$
Let us introduce a permutation matrix $\nu^\ast \in SL(k,\bZ)$
whose $j-$th column equals to
$$^t(\overbrace{0, \cdots,0}^{\nu(j)-1} 
\rlap{\ ${}^{\nu(j)\atop{\hbox{${}^{\vee}$}}}$},1,\;
\overbrace{0, \cdots,0}^{n-\nu(j)}).$$
The last expression $(3.9)$ turns out to be $-\nu^\ast$.
This fact can be seen by the relation
$ \sf L \cdot \sf L^{-1}=id_L$ 
and the Proposition ~\ref{prop21}, 1), $\bf a, b, c$.

We can define a set of indices $^TI_\Lambda$ for $(1.4)^T$
analogous to  $I_\Lambda.$
We have the indices of $I_\Lambda$, $1=j_1 <\cdots <j_n =L-4$
and those of   $^TI_\Lambda,$ $1=i_1 <\cdots <i_n =L-4.$
In this situation we consider two conditions on the CI 
$(1.1),$
$$
\left (\begin {array}{c}
^T\xi_{i_1}(z)\\
\vdots\\
^T\xi_{i_n}(z)\\
\end {array}\right )
= G \cdot V^\Lambda \cdot \left (\begin {array}{c}
^T\xi^{(1)}(z)\\
\vdots\\
^T \xi^{(k)}(z)\\
\end {array}\right ), \leqno(3.10)$$
for linear functions
$(^T\xi^{(1)}(z), \cdots,
^T\xi^{(k)}(z))$ and
$$
\left (\begin {array}{c}
\xi_{j_1}(z)\\
\vdots\\
\xi_{j_n}(z)\\
\end {array}\right )
=^TG \cdot ^TV^\Lambda \cdot \left (\begin {array}{c}
\xi^{(1)}(z)\\
\vdots\\ \xi^{(k)}(z)\\
\end {array}\right ), \leqno(3.10)^T$$
for possibly another $k-$tuple of
linear functions $(\xi^{(1)}(z), \cdots,\xi^{(k)}(z)).$

Due to the condition $(1.11)$
on the set of indices $^TI_\Lambda$, we have for some
permutation matrices $\rho, ^T\rho \in SL(n, \bZ)$,
$$ \rho \cdot V^\Lambda= G^{-1}\cdot \left (^t(\vec g^{(1)}),\cdots,
^t(\vec g^{(k)})\right ), $$
$$ ^T\rho \cdot ^TV^\Lambda= ^TG^{-1}\cdot
\left(^t(^T\vec g^{(1)}),\cdots, ^t(^T \vec g^{(k)})\right ). $$
On these $\rho$ and $^T\rho,$
we impose the following conditions,
$$  G\cdot\rho= ^t (G\cdot\rho), \leqno(3.11)$$
$$  ^TG\cdot ^T\rho=^t(^TG\cdot ^T\rho).\leqno(3.11)^T$$
\par

Under these conditions we formulate the following Theorem that verifies an
hypothesis proposed by \cite{Ber}, \S 3.3 under certain conditions.
\begin{thm}
The Mellin transform $M_{0,\gamma}^0(z)$
of the period integral
$I_{x^0,\gamma}^0(s)$ for the Calabi-Yau CI $(1.1)$
has the following form up to a $\Delta-$ periodic function
in the sense of Definition ~\ref{dfn21},
if it satisfies the conditions
$(3.10),$$(3.10)^T,$$(3.11),$
$(3.11)^T,$
$$ M_{0,\gamma}^0(z(\xi))= \frac{\prod_{\nu=1}^k \prod_{j=1}^{\tilde
\tau_\nu}
\Gamma(\;^Tg_j^{(\nu)}\xi^{(\nu)})}{
\prod_{q=1}^k \Gamma(\sum_{\nu=1}^k
\;^TQ_q^{(\nu)} \xi^{(\nu)})}.
\leqno(3.12)$$ Here
$z(\xi) =(z_1(\xi), \cdots , z_k(\xi))$  is a $k-$tuple
of linear functions in variables $\xi =(\xi^{(1)}, \cdots , \xi^{(k)})$
defined by the relation $(3.10)^T.$

In a symmetric way the Mellin transform
$ M_{0,\;^T\gamma}^0(z(^T\xi))$
for the Calabi-Yau CI, $(1.1)^T$ admits an expression
as follows up to a
$\Delta-$ periodic function,
$$
M_{0,\;^T\gamma}^0(z(\xi))=
\frac{\prod_{\nu=1}^k \prod_{j=1}^{\tau_\nu}
\Gamma(g_j^{(\nu)}\;^T\xi^{(\nu)})}
{\prod_{q=1}^k \Gamma(\sum_{\nu=1}^k
Q_q^{(\nu)} \;^T\xi^{(\nu)})}. \leqno(3.12)^T$$
The functions $^T\xi^{(\nu)}$ are defined by the relation
$(3.10).$
\label{thm31} \end{thm}

To prove the Theorem we prepare the following lemma.

\begin{lem}
Under the conditions imposed on $(1.1)$ in  the
theorem ~\ref{thm31},
the Mellin transform $M_{0,\gamma}^0(z)$
of the period integral $I_{x^0,\gamma}^0(s)$ for the CI  $(1.1)$
admits  up to a $\Delta-$periodic function (in the sense of 
Definition ~\ref{dfn21}) an
expression as follows,
$$M_{0,\gamma}^{0}(z)  = \prod_{i=1}^k \Gamma(z_i)
\prod_{j \in I_\Lambda} \Gamma(\xi_j(z)).$$
\label{lem43}
\end{lem}
{\bf Proof}
It is necessary to show that special solutions
of the system $(1.5)$ satisfy,
$$ \xi_{a^{\nu}-2}(z)
=\xi_{a^\nu-1}(z)=z_{\nu},\;\;\xi_{a^\nu}(z) = 1-z_\nu. \leqno(3.13)$$

To see this, we remark that the system below is a direct consequence of
the relation
${\sf L}{\sf L}^{-1}= id_L,$
$$ \sum_{i=1}^n v_{q,i}^{(\nu)} w_i^{(a^\nu)} +
q_{2\mu-1}^{(a^\nu)} =0, \;\; a^{\mu-1}+1 \leq q \leq a^{\mu-1}+\tau_\mu,
\;\; \mu\in [1,k].\leqno(3.14)$$
$$- \sum_{i=1}^n w_i^{(a^\nu)}= \sum_{i=1}^n w_i^{(a^{\nu}-1)}= 1.$$
The last equality can be deduced from $(3.4)$ and  $(3.5).$
Let  $\mu(j) \in [1,k]$
be an index such that $j$ belongs to $I^{(\mu(j))}.$ 
Then we have
$$ \sum_{i \in I_\Lambda} v_{i,j}^{(\mu(j))} p_q^{(i)} +
p_q^{(a^{\mu(j)-1})} =0,\leqno(3.15)$$
for $q \in [1,k].$

We get the following relation,
$$q_{2\nu-1}^{(a^\nu)}=q_{2\nu}^{(a^\nu)}=1,
\leqno(3.16)$$
$$q_{2\nu-1}^{(a^{\nu'})}=q_{2\nu}^{(a^{\nu'})}=0,\;\;
\rm{for}\; \nu' \not = \nu.
\leqno(3.17)$$
The equality $q_{2\nu}^{(a^{\nu'})}=0$ 
can be derived from the fact that the product of 
the $a^{\nu-1}-$th column of  ${\sf L}^{-1}$ with the
$a^{\nu'}-$th row of ${\sf L}$ is equal to  $0.$
On the other hand, the product of the
$(n+2\nu)-$th column of  ${\sf L}$ with the
$a^{\nu'}-$th row of  ${\sf L}^{-1}$ is equal to 
$q_{2\nu-1}^{(a^{\nu'})}+q_{2\nu}^{(a^{\nu'})}=0.$
This proves  $(3.17).$
The product of the
$a^{\nu}-$th column of  ${\sf L}^{-1}$ with the 
$a^{\nu'}-$th row of ${\sf L}$ is equal to $
q_{2\nu-1}^{(a^{\nu'})} +p_{\nu}^{(a^{\nu'})}=0.$
One deduces the equality  $p_{\nu}^{(a^{\nu'})}=0$
for  $\nu \not = \nu'.$
The product of 
$(a^{\nu}-1)-$th column of ${\sf L}^{-1}$ with the
$(a^{\nu'}-2)-$th row of  ${\sf L}$ is equal to  $
q_{2\nu'-1}^{(a^{\nu}-1)} +p_{\nu'}^{(a^{\nu}-1)}=0.$
As a consequence, we have
$$ p_{\nu'}^{(a^{\nu}-1)}=0 \;\;\rm{for}\;\; \nu \not = \nu', \leqno(3.18)$$
$$ p_{\nu}^{(a^{\nu}-1)}=-q_{2\nu-1}^{(a^{\nu}-1)}  =1,\;\;
p_{\nu}^{(a^{\nu})}=-1,$$
because $ p_{\nu}^{(a^{\nu})}+p_{\nu}^{(a^{\nu}-1)} =0. $
Additionally we see that
$$ q_{r}^{(a^{\nu} -1)}=0,\;\; r \not = 2\nu-1. \leqno(3.19)$$

In summary, by virtue of   $(3.14), (3.16), (3.18),$
$$ \xi_{a^\nu}(z)= \sum_{i=1}^n w_{i}^{(a^{\nu})} + \sum_{r=1}^{2k}
q_{r}^{(a^{\nu})} + \sum_{q=1}^k  p_{q}^{(a^{\nu})} z_q
= -1 +1 + 1 -z_\nu =  1 -z_\nu.$$
On the other hand $(3.14), (3.19), (3.18),$
yield
$$  \xi_{a^\nu-1}(z)= 1+0-1 + z_\nu= z_\nu.$$

As for the function $  \xi_{a^\nu-2}(z)$ it is easy to see that 
all elements of the $ (a^\nu-2)-$th column of  ${\sf L}^{-1}$
consist of zeros except the $(n+2k+\nu)-$th element which is equal to 
$1.$

We thus have the equality,
$$M_{0,\gamma}^{0}(z)  = \prod_{i=1}^k\Gamma(z_i)^2 \Gamma(1-z_i)
\prod_{j \in I_\Lambda} \Gamma(\xi_j(z)) = \prod_{i=1}^k \frac{\pi
}{sin\;\pi z_i} \Gamma(z_i) \prod_{j \in I_\Lambda}
\Gamma(\xi_j(z)).$$ {\bf Q.E.D.}

{\bf Proof of the Theorem ~\ref{thm31}}
Our main task is to show the following relation,
$$
\nu^\ast \cdot
\left (\begin {array}{c}
1-z_1\\
\vdots\\
1-z_k\\
\end {array}\right )
=^T{\hat Q}\left (\begin {array}{c}
\xi^{(1)}(z)\\
\vdots\\ \xi^{(k)}(z)\\
\end {array}\right ), \leqno(3.20)$$
for the permutation matrix $\nu^\ast \in SL(k,\bZ)$
introduced just after the formula $(3.9)$.
To do this, first of all we modify the relation,
$$\;^T{\hat Q}
= ^t(^TV^\Lambda) \cdot
\left (^t(^T\vec g^{(1)}),\cdots,^t(^T \vec g^{(k)})\right )
=\left (\begin {array}{c}
\;^T\vec g^{(1)}\\
\vdots\\ \;^T\vec g^{(k)}\\
\end {array} \right )\cdot^TG^{-1}
\cdot ^t(^T \rho)^{-1}  \cdot
\left(^t(^T\vec g^{(1)}),\cdots, ^t(^T \vec g^{(k)})\right),
\leqno(3.21)$$
that can be derived from $(1.8)^T.$ 
Here we remark that after the definition of
$^T\rho \in SL(n, 
\bZ)$ intrduced just before $(3.11)$ the following relation holds,
$$ ^TV^\Lambda=^T\rho^{-1} \cdot^TG^{-1}\cdot
\left(^t(^T\vec g^{(1)}),\cdots, ^t(^T \vec g^{(k)})\right ). \leqno(3.22) $$
From this relation and $(3.11)^T$ we see that 
$$\left (\begin {array}{c}
\;^T\vec g^{(1)}\\
\vdots\\ \;^T\vec g^{(k)}\\
\end {array} \right )
\cdot ^T \rho^{-1} \cdot^TG^{-1} 
\left(^t(^T\vec g^{(1)}),\cdots, ^t(^T \vec g^{(k)})\right)
= \left (\begin {array}{c}
\;^T\vec g^{(1)}\\
\vdots\\ \;^T\vec g^{(k)}\\
\end {array} \right )\cdot^TV^\Lambda. \leqno(3.23)$$
By virtue of the condition $(3.10)^T$, the following equality holds
$$\left (\begin {array}{c}
\;^T\vec g^{(1)}\\
\vdots\\ \;^T\vec g^{(k)}\\
\end {array} \right )\cdot^TV^\Lambda
\cdot
\left(\begin {array}{c}
\;\xi^{(1)}(z)\\
\vdots\\ \;\xi^{(k)}(z)\\
\end {array} \right )
=
\left (\begin {array}{c}
\;^T\vec g^{(1)}\\
\vdots\\ \;^T\vec g^{(k)}\\
\end {array} \right )\cdot ^TG^{-1}
\left(\begin {array}{c}
\;\xi_{j_1}(z)\\
\vdots\\ 
\;\xi_{j_n}(z)\\
\end {array} \right ).  \leqno(3.24)$$
The combination of $(3.6)$ and $(3.8)$ entails
that the expression $(3.24)$ is equal to
$$
-\left (\begin {array}{c}
\;^T\vec g^{(1)}\\
\vdots\\ \;^T\vec g^{(k)}\\
\end {array} \right )\cdot
^TG^{-1}
\cdot ^t({\sf L}_\Lambda)^{-1} 
\;V^\Lambda \cdot
\left (\begin {array}{c}
1-z_1\\
\vdots\\
1-z_k\\
\end {array}\right )
=
-\left (\begin {array}{c}
\;^T\vec g^{(1)}\\
\vdots\\ \;^T\vec g^{(k)}\\
\end {array} \right )\cdot ^TG^{-1}
\cdot (^T{\sf L}_\Lambda)^{-1}
\;^TV^\Lambda \cdot
\left (\begin {array}{c}
1-z_1\\
\vdots\\
1-z_k\\
\end {array}\right ).$$
From $(3.9)$ it follows that the last expression equals
to $$ 
\nu^\ast 
\cdot
\left (\begin {array}{c}
1-z_1\\
\vdots\\
1-z_k\\
\end {array}\right ).$$
This means $(3.20)$ and consequently,
$$\sum_{\nu=1}^k\;^TQ^{(\nu)}_q \xi^{(\nu)}(z) = 1-z_q,\;\; 1 \leq q \leq k.$$
From the last equality
we can directly derive the relation to be proved 
$$
{\small \prod_{q=1}^k \Gamma(z_q) =
\prod_{q=1}^k
\Gamma(1-\sum_{\nu=1}^k\;^TQ^{(\nu)}_q \xi^{(\nu)}(z))
=\frac{\pi}
{\prod_{q=1}^k \Gamma ( \sum_{\nu=1}^k\;^TQ^{(\nu)}_q \xi^{(\nu)}(z) )
sin \left (\;\pi \sum_{\nu=1}^k\;^TQ^{(\nu)}_q \xi^{(\nu)}(z) \right)}.}$$

On the other hand the lemma ~\ref{lem43}
and the condition $(3.10)^T$ gives us
$$ \xi_{a^{\nu-1}+j}=
\;^Tg_j^{(\nu)}\xi^{(\nu)},\;\; j\in [1, \tau_\nu]$$
that means
$$ \prod_{j \in I_\Lambda} \Gamma(\xi_j(z))
= \prod_{\nu=1}^k \prod_{j =1}^{\tau_\nu}
\Gamma(\;^Tg_j^{(\nu)}\xi^{(\nu)}),$$
which proves $(3.12)$.
The formula $(3.12)^T$ can be proven in a parallel way.
{\bf Q.E.D.}


{
\center{\section{ Duality between monodromy data and Poincar\'e
polynomials}}
}

In connection with the mirror symmetry, we consider the structural algebra
of the CI  $(1.1)$ of dimension $n-k$ denoted by  $X$ $= X_\b1$,
$$ A_X:= \frac{\bC [x]}{(f_1+f_2-1, \cdots f_{2k-1}+f_{2k}-1) \bC
[x]},$$
and a natural filtration on it,
$$ A^j_X := \bigoplus \bC
\{x^\alpha \in A_X; \langle \alpha, \vec g^{(1)}
\rangle =j_1, \cdots,  \langle \alpha, \vec g^{(k)}
\rangle =j_k
 \},$$
with the Poincar\'e polynomial,
$$ P_{A_X}(\lambda) = \sum_{j \in \bZ_{\geq 0}} dim (A_X^j) \lambda_1^{j_1}
\cdots \lambda_k^{j_k}.$$
In an analogous way,
we define corresponding notions of the CI
$Y$ defined by $(1.1)^T,$
$$ A_Y:= \frac{\bC [x]}{(^Tf_1+ ^Tf_2-1, \cdots ^Tf_{2k-1}
+^Tf_{2k}-1) \bC
[x]},$$
$$ A^j_Y :=  \bigoplus\{x^\alpha
\in A_Y; \langle \alpha, \;^T\vec g^{(1)}
\rangle =j_1, \cdots \langle, \alpha, \;^T\vec g^{(k)}
\rangle =j_k \},$$
$$ P_{A_Y}(\lambda) = \sum_{j \in \bZ_{\geq 0}}
dim (A_Y^j) \lambda^j.$$
In this situation the classical result due to
\cite{Dolg} gives us,
$$ P_{A_X}(\lambda) = \frac
{\prod_{q=1}^k \prod_{\nu=1}^k (1-\lambda_\nu^{Q^{(\nu)}_q})}
{\prod_{\nu=1}^k \prod_{j=1}^{\tau_\nu} (1-\lambda_\nu^{g^{(\nu)}_j})},
P_{A_Y}(\lambda) = \frac{\prod_{q=1}^k \prod_{\nu=1}^k
(1-\lambda_\nu^{^TQ^{(\nu)}_q})}{\prod_{\nu=1}^k
\prod_{j=1}^{\tilde\tau_\nu} (1-\lambda_\nu^{^Tg^{(\nu)}_j})}.
\leqno(4.1)$$
Further we introduce the variables
$(t_1, \cdots, t_k) \in \bT^k$
such that
$$ \prod^k_{\nu=1} s_\nu^{^TQ^{(q)}_\nu} = t_q, \;\; 1 \leq q \leq k.$$
If we assume that $rank\;^T{\hat Q} =k$,
this equation is always solvable with respect to the
variables $s=s(t).$
We consider the Mellin inverse transform of
$M_{0,\gamma}^0(z(\xi))$ associated to the CI  (1.1),
$$ U_\alpha(s)=
\int_{\check \Pi_\alpha}
\frac{\prod_{\nu=1}^k \prod_{j=1}^{\tau_\nu}
\Gamma(\;^Tg_j^{(\nu)}\xi^{(\nu)}(z))}{
\prod_{\nu=1}^k \Gamma(\sum_{q=1}^k
\;^TQ_q^{(\nu)} \xi^{(\nu)}(z))} s^{-\bz +\b1} d\bz,$$
where ${\check \Pi}_\alpha \subset \bT^k$  is a cycle avoiding
the singular loci of the integrand.
It is easy to check
$$ U_\alpha(s(t))= det (^T \hat Q)^{-1}\int_{^T \hat Q_\ast^{-1}(\check \Pi_\alpha)}
\frac{\prod_{\nu=1}^k \prod_{j=1}^{\tilde \tau_\nu}
\Gamma(\;^Tg_j^{(\nu)}\xi^{(\nu)})}{
\prod_{\nu=1}^k \Gamma(\sum_{q=1}^k
\;^TQ_q^{(\nu)} \xi^{(\nu)})} t_1^{-\xi^{(1)}}\cdots t_k^{-\xi^{(k)}}
d\xi^{(1)} \wedge \cdots \wedge d\xi^{(k)}. $$
The following system of
differential equations annihilates
the inverse Mellin transform $ U_\alpha(s(t))$,
$$  L_\nu (t_\nu, \vartheta_t)U_\alpha(s(t))=0
, \;\; 1 \leq \nu \leq k,$$
where
$$ L_\nu (t_\nu, \vartheta_t)=
\left(\prod_{j=1}^{\tilde \tau_\nu}\prod_{r=0}^{^Tg_j^{(\nu)}-1}
( - ^Tg_j^{(\nu)} \vartheta_{t_\nu} +r) - t_\nu
\prod_{q=1}^k \prod_{r=0}^{\;^TQ_q^{(\nu)}-1}(\sum_{\mu=1}^k
\;^TQ_q^{(\mu)} \vartheta_{t_\mu} -r) \right), \;\; 1 \leq \nu \leq k.
\leqno(4.2)$$
 We denote by
$\chi_\nu$
the degree of the operator $L_\nu (t, \vartheta_t):$
$\chi_\nu=\sum_{q=1}^k\;^TQ_q^{(\nu)}=
\sum_{j=1}^{\tilde \tau_\nu}\;^Tg_j^{(\nu)}=
\;^Tg_{{\tilde \tau_\nu}+1}^{(\nu)} $
that has already been introduced in $(1.3)^T.$

Here we define the restriction of the operator
$ L_{\nu}(t, \vartheta_t)$ onto the torus
$\bT = \{t \in \bC^k; t_i=0, i \not = \nu\} \setminus \{
t_\nu=0\}$ as follows,
$$ {\tilde L}_{\nu}(t_\nu, \vartheta_{t_{\nu}})
:= \left( \prod_{j=1}^{\tilde \tau_\nu}\prod_{r=0}^{^Tg_j^{(\nu)}-1}
( - ^Tg_j^{(\nu)} \vartheta_{t_\nu} +r) - t_\nu
\prod_{q=1}^k \prod_{r=0}^{\;^TQ_q^{(\nu)}-1}(
\;^TQ_q^{(\nu)} \vartheta_{t_\nu} -r) \right). \leqno(4.2)'$$

On the $\chi_q$ dimensional solution
space of the operator
${\tilde L}_q(t_q, \vartheta_{t_q})$,
we consider the monodromy
$M_q^{(0)} \in GL(\chi_q, \bC) $
(resp.
$M_q^{(\infty)} \in GL(\chi_q, \bC) $)
around the point $t_q=0$ (resp. $t_q=\infty$).
Then we have the following characteristic polynomials of
the monodromy that can be easily calculated from the expression
${\tilde L}_{\nu}
(t_\nu, \vartheta_{t_{\nu}})$,
$$ det (M_q^{(\infty)} - \lambda_q \cdot id_{\chi_q})
= \prod_{\nu=1}^k
(1-\lambda_q^{^TQ^{(q)}_\nu}),
det (M_q^{(0)} - \lambda_q \cdot id_{\chi_q})
= \prod_{j=1}^{\tau_q}
(1-\lambda_q^{^Tg^{(q)}_j}).$$
As a consequence the rational function defined by
$$ M_X(\lambda_1, \cdots \lambda_k): =
\prod^k_{q=1} \frac
{det (M_q^{(\infty)} - \lambda_q \cdot id_{\chi_q})}
{det (M_q^{(0)} - \lambda_q \cdot id_{\chi_q})}$$
has a form
$$M_X(\lambda_1, \cdots \lambda_k)=
\frac{\prod_{q=1}^k \prod_{\nu=1}^k
(1-\lambda_q^{^TQ^{(q)}_\nu})}{\prod_{q=1}^k
\prod_{j=1}^{\tilde \tau_q} (1-\lambda_q^{^Tg^{(q)}_j})}. \leqno(4.3)$$
For the rational function
$M_Y(\lambda_1, \cdots \lambda_k)$
defined in a parallel way to the function
$M_X(\lambda_1, \cdots \lambda_k),$ we have
$$M_Y(\lambda_1, \cdots \lambda_k)=
\frac{\prod_{q=1}^k \prod_{\nu=1}^k
(1-\lambda_q^{Q^{(q)}_\nu})}{\prod_{q=1}^k
\prod_{j=1}^{\tau_q} (1-\lambda_q^{g^{(q)}_j})}.$$

Let $\bar Y$ be the compactification of CI $Y \subset \bT^{n}$
in the product of quasihomogeneous projective spaces
$\bP:=\bP^{(\tau_1)}_{(^Tg_1^{(1)},\ldots,
^Tg_{\tau_1+1}^{(1)})} \times
\ldots \times \bP^{(\tau_k)}_{(^Tg_1^{(k)},\ldots,
^Tg_{\tau_k+1}^{(k)})}.$
The coherent sheaf on $\bP$
${\mathcal O}_{\bP}(\zeta),$
$\zeta=(\zeta_1, \cdots, \zeta_k)$
is defined by sections on the open set
$U_I=\{x \in \bC^n; x_i \not=0, i \in I   \}$
that are given by
$$ \Gamma(U_I, {\mathcal O}_{\bP}(\zeta)):= \oplus
\bC \{x^\alpha; \alpha=(\alpha_1,
\cdots,\alpha_n) \in \bZ^n,
\alpha_i \geq 0 \;for \; i \not \in I,\langle ^Tg^{(q)}, \alpha \rangle =
\zeta_q, 1 \leq q \leq k \}.$$

We define the coherent sheaf
${\mathcal O}_{\bar Y}(\zeta) $
by sections on the open set
$U_I \cap Y,$
$$ {\small \Gamma(U_I, {\mathcal O}_{\bar Y}(\zeta)):=
\oplus
\bC \{x^\alpha \in A_Y; \alpha= (\alpha_1,
\cdots,\alpha_n) \in \bZ^n,
\alpha_i \geq 0 \;\;for \;\;i \not \in I, \langle ^Tg^{(q)},
\alpha \rangle =
\zeta_q, 1 \leq q \leq k \}.}$$
We introduce the Euler characteristic for this sheaf,
$$ \chi({\mathcal O}_{\bar Y}(\zeta)):=
\sum^{n-k}_{i=0}(-1)^i dim H^i({\mathcal O}_{\bar Y}(\zeta)).$$

After \cite{Dolg}, \cite{Gol} the Poincar\'e
polynomial of the
Euler characteristics,
$$ P{\mathcal O}_{\bar Y}(t_1, \cdots, t_k)
:= \sum_{\zeta \in (\bZ_{\geq 0})^k}\chi({\mathcal O}_{\bar Y}(\zeta))
t_1^{\zeta_1} \cdots t_k^{\zeta_k}$$
admits an expression as follows,
$$ P{\mathcal O}_{\bar Y}(t_1, \cdots, t_k)=
\frac{\prod_{q=1}^k \prod_{\nu=1}^k
(1-t_q^{^TQ^{(q)}_\nu})}{\prod_{q=1}^k
\prod_{j=1}^{\tilde \tau_\nu} (1-t_q^{^Tg^{(q)}_j})}.
\leqno(4.4)$$
For the sheaf
${\mathcal O}_{\bar X}(\zeta)$
defined analogously to
${\mathcal O}_{\bar Y}(\zeta)$, we consider
the Poincar\'e polynomial of the
Euler characteristics
$$ P{\mathcal O}_{\bar X}(t_1, \cdots, t_k)
:= \sum_{\zeta \in (\bZ_{\geq 0})^k}\chi({\mathcal O}_{\bar X}(\zeta))
t_1^{\zeta_1} \cdots t_k^{\zeta_k}.$$
We have the following expression after \cite{Dolg} and \cite{Gol},
$$ P{\mathcal O}_{\bar X}(t_1, \cdots, t_k)=
\frac{\prod_{q=1}^k \prod_{\nu=1}^k
(1-t_q^{Q^{(q)}_\nu})}{\prod_{q=1}^k
\prod_{j=1}^{\tau_\nu} (1-t_q^{g^{(q)}_j})}.$$
If we compare $(4.1),$ $(4.3)$
and  $ (4.4)$
we get the following statement.
\begin{thm}
For the Calabi-Yau CI's  $X$ and $Y$ defined by  $(1.1)$ and $(1.1)^T$
we have the following relations,
$$M_X(\lambda_1,\cdots, \lambda_k)=  P{\mathcal O}_{\bar Y}
(\lambda_1,\cdots,\lambda_k)
= P_{A_Y}(\lambda_1,\cdots,\lambda_k),$$
$$M_Y(\lambda_1,\cdots, \lambda_k)=
P{\mathcal O}_{\bar X}(\lambda_1,\cdots,
\lambda_k)
= P_{A_X}(\lambda_1,\cdots,\lambda_k),$$
if they satisfy sufficient conditions of Theorem ~\ref{thm31}.
\label{thm432}
\end{thm}

In this situation
one can easily derive the existence of the following symmetry
from the Theorem~\ref{thm31}.
\begin{thm}
Assume that  $X$, $(1.1)$ and $Y$, $(1.1)^T$ satisfy all conditions
imposed on them in Theorem ~\ref{thm31}, then
there is a symmetry between geometric symmetry and
quantum symmetry
as follows.
$${\mathcal Q}_{X} \cong \prod^{k}_{q=1}{\bZ_{\bf {\bar Q}^{(q)}}}
\cong {\mathcal G}_{ Y},$$
$${\mathcal Q}_{ Y} \cong \prod^{k}_{q=1}{\bZ_{\bf^T{\bar Q}^{(q)}}}
\cong {\mathcal G}_{X}.$$
\label{thm40}
\end{thm}
{\bf Proof} The isomorphism
${\mathcal Q}_{X} \cong \prod^{k}_{q=1}{\bZ_{\bf{\bar Q}^{(q)}}}$
and its analogy on $Y_s$ is clear from the quasihomoeneiety of the
systems $(1.1)$ and $(1.1)^T$.
The existence of the cyclic action
$\prod^{k}_{q=1}{\bZ_{\bf^T{\bar Q}^{(q)}}}$ on $X_s$
can be read off from the monodromy actions of $M_q^{(\infty)}, 1 \leq q \leq k$
on the solutions to the system $(4.2)$
that are defined along vanishing cycles on $X_s.$
Here we remark that the monodromy group is a subgroup of $Aut(X_s)$
in view of the integer power transform from $(s_1, \cdots, s_k)$
to $(t_1, \cdots, t_k)$ defined just after $(4.1)$.
We remark here that the monodromy actions of $M_q^{(0)}, 1 \leq q \leq k$
consist part of the cyclic action
$\prod^{k}_{q=1}{\bZ_{\bf ^T{\bar Q}^{(q)}}}$
because every $^Tg^{(q)}_j,$ $j=1, \cdots, \tilde \tau_q$
divides $\bf ^T{\bar Q}^{(q)}$.
As for the monodromy of the solution to the system $(4.2)'$
an analogous argument holds.
{\bf Q.E.D.}

\begin{remark}
{\em As we have not calculated the global monodromy of the integrals
$U_\alpha(s(t))$, $(4.3)$, it is not proper to talk about
the existence of a mirror symmetry between
$X$ and $Y.$ We gathered, however, the monodromy
data which correspond to  certain
limit of the values
$t_\ell \rightarrow 0, \ell \not = \nu.$
Roughly speaking, this procedure can be interpreted as the
selection of a special face of the Newton polyhedron $\Delta(F(x,s,y))$
(1.4) for the calculus of the monodromy data.

In the article \cite{CKYZ},
authors studied the period integral of 
$X$ at certain limit of the parameter values
and they deduced informations on the analogies 
Gromov-Witten invariants of  $Y$. 
They call this duality ``local mirror symmetry''.
It is probable that our Theorem
~\ref{thm432} is one of numerous 
aspects of the local mirror symmetry.}
\end{remark}

{
\center{\section{ Dual nef-partition interpretation}}
}

In this section we show that under certain condition the transposition 
mirror construction corresponds to the notion of the dual nef-partion due to
L.Borisov \cite{BB}.

First of all we consider the following set of vectors defined after 
$(1.1)$, $(1.2).$
$$ \leqno(5.1)
\begin {array}{ccc}
\vec v^{(q)}_1-\vec v^{(q)}_{\tau_q+1}&&\\
\vdots& &1 \leq q \leq k\\
\vec v^{(q)}_{\tau_q}-\vec v^{(q)}_{\tau_q+1}&&\\
\end {array}$$
This gives rise to a partition of $I_\Lambda$ with
$ \sharp I_\Lambda =n.$
We set 
$$\Delta_q := {\rm convex\;hull\;of}\; (\{0\} \cup \bigcup_{j=1}^{\tau_q} \{ 
\vec v^{(q)}_j-\vec v^{(q)}_{\tau_q+1} \} ). \leqno (5.2)$$
Further on we impose the following condition on 
the Minkowski sum of $\Delta_q$,
$$ dim(\Delta_1 + \cdots +\Delta_k)= n-k. \leqno(5.3)$$
We introduce the
$(n-k)$ dimensional integral lattice
$$ V_\bZ = \{ x  \in \bZ^n; \langle x, \vec g^{(q)} \rangle
=0, 1 \leq q \leq k\}.$$
After  this notation each $\Delta_q$ is located on
$V_\bR = V_\bZ \times \bR \cong \bR^{n-k}.$
Consequently $\Delta_1 + \cdots +\Delta_k \subset V_\bR.$
There exists a piecewise linear function $\phi_q$ such that
$ \phi_q(y_1+y_2) \leq  \phi_q(y_1) +  \phi_q(y_2).$
Namely we shall define it as follows,
$$  \phi_q(y) = - min_{ x \in \Delta_q} \langle x,y\rangle. \leqno(5.4)$$
We construct a set of polyhedra dual to the set 
$\{\Delta_1, \cdots , \Delta_k  \},$
$$ \Delta_q^\ast = \{\vec m \in \bR^n; (5.5)_\ell, (5.5)'_\ell
\}, \;\; 1 \leq \ell \leq k. \leqno(5.5)_0$$
The conditions $(5.5)_\ell, (5.5)'_\ell$ look like the following,
$$  \langle\vec m,\vec v^{(\ell)}_j-\vec v^{(\ell)}_{\tau_\ell+1}\rangle \geq -1, 
1 \leq j \leq  \tau_\ell. \leqno{(5.5)_\ell}$$
$$  \langle\vec m,\vec v^{(q)}_i-\vec v^{(q)}_{\tau_q+1}\rangle \geq 0, 1 \leq i \leq  \tau_q,
q \not = \ell. \leqno{(5.5)_\ell}$$
We can construct the set 
$ \{\Delta_1^\ast, \cdots , \Delta_k^\ast  \}  $
by means of the vertices vectors $\vec m^{(\ell)}_1, \cdots, \vec m^{(\ell)}
_{\tau_\ell}, 1 \leq \ell \leq k$ which are defined by the following 
set 
of equalities and inequalities.
$$ \langle\vec m^{(\ell)}_r,\vec v^{(\ell)}_j-\vec v^{(\ell)}_{\tau_\ell+1}\rangle = -1, 
j \not = r, 1 \leq j \leq  \tilde \tau_\ell. \leqno(5.6)_1$$ 
$$ \langle\vec m^{(\ell)}_r,\vec v^{(\ell)}_r-
\vec v^{(\ell)}_{\tau_\ell+1}\rangle = -1. \leqno(5.6)_2  $$
$$ \langle\vec m^{(\ell)}_r,\vec v^{(q)}_j-\vec v^{(q)}_{\tau_q+1}\rangle =0. \leqno(5.6)_3$$
for $q \in [1,k]$ and $j \in [1,\tau_q] \setminus j_q$ for some index
$j_q$ associated to $\vec m^{(\ell)}_r.$
$$\langle\vec m^{(\ell)}_r,\vec v^{(q)}_{j_q}-\vec v^{(q)}_{\tau_q+1}\rangle \geq 
0. \leqno(5.6)_4$$
By virtue of the condition $(5.3)$
we have a set of $(n-k)$ independent 
equations corresponding to the equalities
$(5.6)_1$,$(5.6)_3$ above.

We see that we may choose as $\vec m^{(\ell)}_r \in M_\bR
=M_\bZ\otimes \bR$ induced from the mapping
$$ pr:\bZ^n \rightarrow \frac{\bZ^n}
{\sum^k_{q=1} \bZ \vec g^{(q)} } := M_\bZ \cong \bZ^{n-k},$$
in view of the relation $(1.2).$

Therefore we may uniquely detemine the set of
vectors $\{\vec m^{(\ell)}_r \} \in M_\bR$ as solutions to a system of
 $(n-k)$ equations $(5.6)_1$,$(5.6)_3$ under certain compatibility condition.

To formulate this compatibility condition, let us  denote by 
$$ {\sf P} = \;^t(\vec m^{(1)}_1, \cdots, \vec m^{(1)}_{\tau_1}, \vec m^{(2)}_1,
\cdots, \vec m^{(k)}_{\tau_k}   ),$$
a $n \times n$ matrix whose $b^{\ell-1}+r$ th column corresponds to  
$m^{(\ell)}_r$.
The system $(5.6)_\ast$ can be realized by the following matrix
equation if it is solvable,
$$ ({\sf L}_\Lambda - \;\rho \cdot V^\Lambda \cdot ^t V^\Lambda)\cdot {\sf P}
= \;^T{\sf L}_\Lambda -\;^T\rho \cdot ^TV^\Lambda \cdot ^t(^TV^\Lambda).
\leqno(5.7)$$ The solvability of this equation can be understood
as the compatibility mentioned above.
Let us formulate a sufficient condition for the solvability of $(5.7)$.

\begin{lem}
Let us assume that all conditions imposed on
$(1.1)$ and $\;^T(1.1)$ in Theorem ~\ref{thm31}
are satisfied. Furthermore we assume that it is possible to make
$\lambda=id_n$ in $(1.12)$ by means of the rearrangements of rows and columns
in $ {\sf L}_\Lambda$. Assume that $G=id_n$ (resp. $^TG=id_n$)
in $(1.9)$ (resp.  $^T(1.9)$). Then the  equation $(5.7)$
is solvable with respect to ${\sf P}$. 
\label{lem52}
\end{lem}
{\bf proof}
The existence of $k$ linearly independent eigenvectors 
$\vec v^{(0)}_{1}, \cdots , \vec v^{(0)}_{k}$
$\in \bR^n$ such that
$$ 
^t\vec v^{(0)}_{\ell}\cdot({\sf L}_\Lambda - 
\;\rho \cdot V^\Lambda \cdot ^t V^\Lambda)
= ^t\vec v^{(0)}_{\ell}\cdot 
(\;^T{\sf L}_\Lambda -\;^T\rho \cdot ^TV^\Lambda \cdot ^t(^TV^\Lambda))
=\vec 0 \in \bR^k,$$
is a necessary condition for the solvability of the  equation $(5.7)$.
It is also a sufficient condition for the solvability as
the following relations show, 
$$^t({\sf L}_\Lambda - 
\;\rho \cdot V^\Lambda \cdot ^t V^\Lambda) \in  GL(W),$$
$$\;^t{\sf P}\cdot^t({\sf L}_\Lambda - 
\;\rho \cdot V^\Lambda \cdot ^t V^\Lambda)= 
\;^t (\;^T{\sf L}_\Lambda -\;^T\rho \cdot ^TV^\Lambda \cdot ^t(^TV^\Lambda))
\in GL(W),$$
for the vector space $W := \frac{\bR^n}{\sum^k_{\ell=1} \bR
\vec v^{(0)}_{\ell} } \cong \bR^{n-k}.$

We will see further that
$$^t({\sf L}_\Lambda - 
\;\rho \cdot V^\Lambda \cdot ^t V^\Lambda)\cdot V^\Lambda
= ^t (\;^T{\sf L}_\Lambda -\;^T\rho \cdot ^TV^\Lambda \cdot ^t(^TV^\Lambda))
\cdot V^\Lambda =0 \in End(\bR^{n}, \bR^k),$$
under the imposed conditions.

Firt we remark that 
$\;^t(\rho \cdot V^\Lambda) \cdot V^\Lambda = \;^t \hat Q$ in view of $(1.8)$
due to the condition
$G=id_n.$
Thus it is enough to show the equality
$$ (^t{\sf L}_\Lambda)^{-1}V^\Lambda\;^t \hat Q= V^\Lambda. \leqno(5.8)$$
The left hand side of  $(5.8)$
is, in its turn,  equal to,
$$ - (\;^t \vec P_1,\cdots ,\;^t \vec P_k)\cdot\;^t \hat Q,$$
by virtue of $(3.6). $ Let us introduce a matrix,
$$ \tilde {\sf P} = \left (\begin {array}{ccc}
\tilde p^{(1)}_1& \cdots&\tilde p^{(k)}_1 \\
\vdots& \cdots&\vdots\\
\tilde p^{(1)}_k& \cdots&\tilde p^{(k)}_k \\
\end {array}\right ),\leqno(5.9)$$
for $\;^T\xi^{(q)}(z) = \sum^k_{r=1} \tilde p^{(q)}_r(1-z_r).$
With this notation the matrix $(5.9)$ equals to 
$$-V^\Lambda \tilde {\sf P}\;^t \hat Q, \leqno(5.10)$$
if we assume $(3.10)$ for $G=id_n.$
Under the conditions imposed on 
$(1.1)$ and $\;^T(1.1)$ in Theorem ~\ref{thm31} we have
$$
\left (\begin {array}{c}
1-z_1\\
\vdots\\
1-z_k\\
\end {array}\right )
={\hat Q}\left (\begin {array}{c}
\;^T\xi^{(1)}(z)\\
\vdots\\ \;^T\xi^{(k)}(z)\\
\end {array}\right ),$$
which is a mere analogy to $(3.10).$ In making use of this equality we see that
$(5.10)$ is equal to $V^\Lambda.$

To see the equality 
$$^t (\;^T{\sf L}_\Lambda -\;^T\rho \cdot ^TV^\Lambda \cdot ^t(^TV^\Lambda))
\cdot V^\Lambda =0 \in End(\bR^{n}, \bR^k),$$
first we prove the equality 
$$^t(\;^T{\sf L}_\Lambda -\;^T\rho \cdot ^TV^\Lambda \cdot ^t(^TV^\Lambda))
\cdot^TV^\Lambda =0,$$
in a way parallel to the proof of $(5.8).$
Further we see that $^TV^\Lambda= V^\Lambda$
under the condition $\lambda =id_n.$
{\bf Q.E.D.}

As the matrices ${\sf L}_\Lambda - \;\rho \cdot V^\Lambda \cdot ^t V^\Lambda$
and $\;^T{\sf L}_\Lambda -\;^T\rho \cdot ^TV^\Lambda \cdot ^t(^TV^\Lambda)$
have rank $(n-k)$ because of the condition
$(5.3)$, the columns of the matrix $\sf P$ are determined
as elements in $ M_\bR.$ We recall that the matrices 
 $\rho \cdot V^\Lambda$ and
$\;^T\rho \cdot ^TV^\Lambda$ have been introduced to formulate 
the conditions $(3.11)$ and $(3.11)^T.$
It is clear that the columns of the solution ${\sf P}$ to this
equation satisfy the equations $(5.6)_\ast$ and, in particular, determine $j_q$
of $(5.6)_4$.

The vectors $(5.1)$ admit another interpretation. 
We introduce unit vectors
$$\epsilon_q = (0,\cdots, \rlap{\ ${}^{q\atop{\hbox{${}^{\vee}$}}}$}1, 0, \cdots,0)
\in \bZ^k, 1 \leq q \leq k. $$
Let us consider
a $(n+k)$ dimensional cone $\sigma$
with $(n+k)$ generators
$$\tilde v^{(q)}_0 = (0,\cdots,0, \epsilon_q) \in \bZ^{n+k},$$
$$ \tilde v^{(q)}_j = (\vec v^{(q)}_j-\vec v^{(q)}_{\tau_q+1}, \epsilon_q) \in \bZ^{n+k}, 
1 \leq j \leq \tau_q, 1 \leq q \leq k.$$ 
The dual cone $\check \sigma$
to the cone $\sigma=\bR_{\geq 0}\langle\tilde v^{(1)}_0, 
\tilde v^{(1)}_1, \cdots, \tilde v^{(1)}_{\tau_1},
\cdots, \tilde v^{(k)}_{\tau_k}\rangle$ is defined as
$$ \check \sigma = \{y \in \bR^{n+k}; \langle x,y\rangle \geq 0\; {\rm for\;all}
\;x \in \sigma\}.$$
The generators of the dual cone $\check \sigma$
are given by the vectors,
$$ \tilde  m^{(\ell)}_0 := (0,\cdots,0, \epsilon_\ell),$$
$$ \tilde  m^{(\ell)}_r := 
(\vec m^{(\ell)}_r, \epsilon_\ell), 1 \leq r \leq \tau_\ell, 
1 \leq \ell \leq k$$ with $\vec m^{(\ell)}_r$
satisfying the equations $(5.6)_\ast.$ It is evident that
the following inequalities hold for the above vectors,
$$ \langle\tilde v^{(q)}_j, \tilde  m^{(\ell)}_r\rangle \geq 0.$$
Conversely all vectors $\tilde  m$
satisfying the conditions 
$$ \langle\tilde v^{(q)}_j, \tilde  m \rangle \geq 0$$
for every $\tilde v^{(q)}_j$ must be a linear combination
of $\tilde  m^{(\ell)}_r$'s with positive coefficients as
they form a basis of $M_\bR \times\bR^{k}$ the dual space to
$V_\bR \times\bR^{k}.$
Thus we get the generators of the dual cone $\check \sigma$
by means of the vertices of  $\{\Delta_1^\ast, \cdots , \Delta_k^\ast  \}.$
This is an example of realization of \cite{BB}, Theorem 4.6.  

Let us formulate a statement on the dual partition
$ \{\Delta_1^\ast, \cdots , \Delta_k^\ast  \} $
that can be characterized as a dual nef-partition 
to $ \{\Delta_1, \cdots , \Delta_k  \}  $.
We spare the space to formulate the definition of the notion of 
dual nef-partition (see  \cite{BB} Definition 4.2)
by  making it clear in the proof 
of the folowing statement.

\begin{prop}
We assume the following two conditions on $(1.1)$ and $^T(1.1)$.
$\ba$. $\Delta_q \cap \Delta_\ell =\{0\}$ for $q \not = \ell$. 
$\bb$. It is possible to choose
an integer entry matrix $\sf P$ in the matrix equation $(5.7)$.
Then the following statements hold.

1.  $ \{\Delta_1, \cdots , \Delta_k  \}  $ is a nef-partition of the Minkowski
sum $ \Delta_1+ \cdots+ \Delta_k$. 

2. $ \{\Delta_1^\ast, \cdots , \Delta_k^\ast  \}  $
is a nef-partition of the Minkowski sum $ \Delta_1^\ast+ 
\cdots+\Delta_k^\ast  $ dual to  $ \{\Delta_1, \cdots , \Delta_k  \}$.

3.The transposed polynomials $\;^Tf_{2q-1}/(\;^Tf_{2q}-1)$
are obtained from $\Delta_q^\ast$ by means of the torus closed embedding,
$$\begin {array}{ccc}
\prod 
\bP^{\tilde \tau_q-1}_{(\;^Tg^{(q)}_1,\cdots, 
\;^Tg^{(q)}_{\tilde \tau_q})}& \rightarrow &\bP^{n-1}\\
(x_1, \cdots, x_n)& \mapsto &(y_1, \cdots, y_n)\\
& & =(\prod^k_{\ell=1}\prod_{1 \leq j \leq \tau_\ell}
x_{b^{\ell-1}+j}^{\langle \vec v^{(\ell)}_{j}-\vec v^{(\ell)}_{\tau_\ell+1}, \be_1\rangle },
\cdots, \prod^k_{\ell=1}\prod_{1 \leq j \leq \tau_\ell}
x_{b^{\ell-1}+j}^{\langle \vec v^{(\ell)}_{j}-\vec v^{(\ell)}_{\tau_\ell+1}, \be_n\rangle}),
\end {array}
$$
where $\be_i= (0,\cdots, \rlap{\ ${}^{i\atop{\hbox{${}^{\vee}$}}}$}1, 0, \cdots,0)$
$\in \bZ^n.$
\label{prop51}
\end{prop}

{\bf Proof}

1. The condition ${\sf P} \in End(\bZ)$
entails the integral linear property of the function
$ \phi_\ell(y)$ defined in $(5.4).$ The convexity of 
$\phi_\ell : M_\bR \rightarrow \bR$
is guaranteed by the following facts. First the number of vertices in
$\Delta_\ell$ is less than $(n-k+1)= dim M_\bR +1$. Secondly 
the condition $\bf a$.

It is easy to see from the equation $(5.7)$ that
$$ \phi_q( \vec m^{(\ell)}_r) = \delta_{q,\ell}, 1 \leq r \leq \tau_\ell. 
\leqno(5.11)$$
The existence of an integral  upper convex piecewise linear function $\phi_q,$
$1 \leq q \leq k$ satisfying $(5.11)$ is equivalent to  the definition of 
nef- partition $ \{\Delta_1, \cdots , \Delta_k  \}$ of the Minkowski
sum $\Delta_1+ \cdots+ \Delta_k$ whose dimension is equal to $dim V_\bR$
after $(5.3)$. 

2. After the condition $\bb$, 
$$\psi_\ell(x) = - min_{y \in \Delta^\ast} 
\langle x,y\rangle,$$
is an integral piecewise linear function on $V_\bR.$
The relation $\Delta_q^\ast \cap \Delta_\ell^\ast =\{0\}$ for $q \not = \ell$
is clear from the existence of the funciton $(5.11).$
The number of vertices in
$\Delta_\ell^\ast$ is less than $(n-k+1)= dim V_\bR +1$.
Thus we see that $\psi_\ell(x)$ is an upper convex function.
This shows that $ \{\Delta_1^\ast, \cdots , \Delta_k^\ast  \}$ 
is a nef-partition of the Minkowski
sum $ \Delta_1^\ast + \cdots+ \Delta_k^\ast$. The existence of
an integral  upper convex piecewise linear function $(5.11)$
shown before means that $ \{\Delta_1^\ast, \cdots , \Delta_k^\ast  \}$ 
is a dual nef-partition to $ \{\Delta_1, \cdots , \Delta_k  \}.$

3.First of all we remark that the transposition construction entails,
$$ \sum^k_{q=1} \;^T\tilde g^{(q)}_j (v^{(\nu(q))}_j-
\vec v^{(\nu(q))}_{\tilde
\tau_q+1})=0.$$
By virtue of the convexity of the function $\sum_{\ell=1}^k\phi_\ell$
(see \cite{oda} \S 2.3) or equivalently 
$$ \Delta_1 + \cdots+ \Delta_k = \{x \in V_\bR; \langle x,y\rangle \geq -
\sum_{\ell=1}^k\phi_\ell(y), y \in M_\bR \},$$ the mapping
$(5.7)$ is a closed torus embedding. Thus 
a polynomial whose Newton diagram equals to 
$\Delta^\ast_\ell$ in $(y_1, \cdots, y_n)$ 
variables coincides with a polynomial obtained 
as a sum of monomials 
with exponents from the rows of the RHS of $(5.7)$ i.e.
$\prod^k_{\ell=1}\prod_{1 \leq j \leq \tilde \tau_\ell}
x_{b^{\nu(\ell-1)}+j}^{ \;^T\vec v^{(\ell)}_{j}-\;^T
\vec v^{(\ell)}_{\tilde\tau_\ell+1}}$.
Further argument is parallel to that in \cite{BK}, \S 3.
{\bf Q.E.D.}

{
\center{\section{ Examples }}
}

{\bf Example 6.1, Schimmrigk variety}
As a simple, but non-trivial example we recall the 
following example whose period integral has been studied
in \cite{Ber},
$$ f_1(x)= \sum_{i=0}^3 x_i^3,
f_2(x)= x_1 x_2 x_3 +1,$$
$$ f_3(x)= \sum_{i=1}^3 x_i x_{i+3}^3,
 f_4(x)= x_0 x_4 x_5 x_6 +1.$$
We then have the matrices below after the notation $(1.5)$,
$${\sf L}=
\left (\begin {array}{ccccccccccccc} 3&0&0&0&0&0&0&1&0&0&0&0&0
\\\noalign{\medskip}0&3&0&0&0&0&0&1&0&0&0&0&0\\\noalign{\medskip}0&0&3
&0&0&0&0&1&0&0&0&0&0\\\noalign{\medskip}0&0&0&3&0&0&0&1&0&0&0&0&0
\\\noalign{\medskip}0&0&0&0&0&0&0&1&0&0&0&1&0\\\noalign{\medskip}0&1&1
&1&0&0&0&0&1&0&0&0&0\\\noalign{\medskip}0&0&0&0&0&0&0&0&1&0&0&0&0
\\\noalign{\medskip}0&1&0&0&3&0&0&0&0&1&0&0&0\\\noalign{\medskip}0&0&1
&0&0&3&0&0&0&1&0&0&0\\\noalign{\medskip}0&0&0&1&0&0&3&0&0&1&0&0&0
\\\noalign{\medskip}0&0&0&0&0&0&0&0&0&1&0&0&1\\\noalign{\medskip}1&0&0
&0&1&1&1&0&0&0&1&0&0\\\noalign{\medskip}0&0&0&0&0&0&0&0&0&0&1&0&0
\end {array}\right ).$$

$${\small {\sf L}^{-1}=
\left (\begin {array}{ccccccccccccc} 1/3&-1/9&-1/9&-1/9&0&1/3&-1/3&0&0
&0&0&0&0\\\noalign{\medskip}0&2/9&-1/9&-1/9&0&1/3&-1/3&0&0&0&0&0&0
\\\noalign{\medskip}0&-1/9&2/9&-1/9&0&1/3&-1/3&0&0&0&0&0&0
\\\noalign{\medskip}0&-1/9&-1/9&2/9&0&1/3&-1/3&0&0&0&0&0&0
\\\noalign{\medskip}-1/9&-1/27&{\frac {2}{27}}&{\frac {2}{27}}&0&-1/9&
1/9&2/9&-1/9&-1/9&0&1/3&-1/3\\\noalign{\medskip}-1/9&{\frac {2}{27}}&-
1/27&{\frac {2}{27}}&0&-1/9&1/9&-1/9&2/9&-1/9&0&1/3&-1/3
\\\noalign{\medskip}-1/9&{\frac {2}{27}}&{\frac {2}{27}}&-1/27&0&-1/9&
1/9&-1/9&-1/9&2/9&0&1/3&-1/3\\\noalign{\medskip}0&1/3&1/3&1/3&0&-1&1&0
&0&0&0&0&0\\\noalign{\medskip}0&0&0&0&0&0&1&0&0&0&0&0&0
\\\noalign{\medskip}1/3&-1/9&-1/9&-1/9&0&0&0&1/3&1/3&1/3&0&-1&1
\\\noalign{\medskip}0&0&0&0&0&0&0&0&0&0&0&0&1\\\noalign{\medskip}0&-1/
3&-1/3&-1/3&1&1&-1&0&0&0&0&0&0\\\noalign{\medskip}-1/3&1/9&1/9&1/9&0&0
&0&-1/3&-1/3&-1/3&1&1&-1\end {array}\right ).}$$

After the construction
$(1.4)^T$ based on the transposed matrix
$\;^t(\sf L)$ it is easy to see that
$ \;^Tf_\ell(x) = f_\ell(x), 1 \leq \ell \leq 4.$
The matrix $\;^T{\sf L}^{-1}= {\sf L}^{-1}$ 
gives us linear functions,
$$\xi^{(1)} = -\frac{1}{3}(z_1-1) +\frac{1}{9}(z_2-1)
, \xi^{(2)} = -\frac{1}{3}(z_2-1).$$
For $f(x)$ as well as for $^Tf(x)$
we can calculate the Mellin transform of the period integral,
$$ M_{0,\gamma}^{0}(z)=
\Gamma( -\frac{1}{3}(z_1-1) +\frac{1}{9}z_2)^3
\Gamma(-\frac{1}{3}(z_2-1))^4 \Gamma(z_1)\Gamma(z_2)
= \frac{\Gamma(\xi^{(1)})^3\Gamma(\xi^{(2)})^4}{\Gamma(3\xi^{(1)}
+\xi^{(2)}) \Gamma(3\xi^{(2)}) },$$
up to $27-$periodic functions.

Analogously we can look at the CI defined on $\bT^{2n+1}$
$$ f_1(x)= \sum_{i=0}^n x_i^n,
f_2(x)= x_1 x_2 \cdots x_n +1,$$
$$ f_3(x)= \sum_{i=1}^n x_i x_{i+n}^n,
 f_4(x)= x_0 x_{n+1} x_{n+2}\cdots x_{2n} +1,$$
whose period integral can be expressed through its Mellin transform,
$$ M_{0,\gamma}^{0}(z)=
 \frac{\Gamma(\xi^{(1)})^n\Gamma(\xi^{(2)})^{n+1}}{\Gamma(n\xi^{(1)}
+\xi^{(2)}) \Gamma(n\xi^{(2)}) },$$
up to $n^n-$periodic functions.

It is worthy to notice that our matrix
$\sf L$ satisfies an interesting condition below.

{\bf The magic square condition}
For each fixed $q   \in [1,k]$, we can find a single valued mapping
$\sigma: b \in I_\Lambda  \rightarrow \{1,
\cdots,n \}$ such that
$$ p_q^b= w_{\sigma(b)}^{a^q}, {\rm{for\; all}}\;b \in I_\Lambda.$$

This condition plays central r\^ole 
in the interpretation of the strange duality found by 
Arnol'd on the interchange between
Gabrielov number and  Dolgachev number from the point of view
of the mirror symmetry \cite{Kob}.

{\bf Example 6.2}
We consider an example of an hypersurface studied in
\cite{Ber}. We have the following data after the notations above,
$$ f_1(x)= x_1^7 + x_2^7x_4+x_3^7x_5+ x_4^3+x_5^3,\;\;
f_2=x_1x_2x_3x_4x_5.$$

$$
{\sf L}=\left [\begin {array}{cccccccc} 7&0&0&0&0&1&0&0\\\noalign{\medskip}0&7
&0&1&0&1&0&0\\\noalign{\medskip}0&0&7&0&1&1&0&0\\\noalign{\medskip}0&0
&0&3&0&1&0&0\\\noalign{\medskip}0&0&0&0&3&1&0&0\\\noalign{\medskip}0&0
&0&0&0&1&0&1\\\noalign{\medskip}1&1&1&1&1&0&1&0\\\noalign{\medskip}0&0
&0&0&0&0&1&0\end {array}\right ],$$

$$
{\sf L}^{-1}=
\left [\begin {array}{cccccccc} {\frac {6}{49}}&-1/49&-1/49&-{\frac {2
}{49}}&-{\frac {2}{49}}&0&1/7&-1/7\\\noalign{\medskip}-{\frac {2}{147}
}&{\frac {19}{147}}&-{\frac {2}{147}}&-{\frac {11}{147}}&-{\frac {4}{
147}}&0&2/21&-2/21\\\noalign{\medskip}-{\frac {2}{147}}&-{\frac {2}{
147}}&{\frac {19}{147}}&-{\frac {4}{147}}&-{\frac {11}{147}}&0&2/21&-2
/21\\\noalign{\medskip}-1/21&-1/21&-1/21&{\frac {5}{21}}&-2/21&0&1/3&-
1/3\\\noalign{\medskip}-1/21&-1/21&-1/21&-2/21&{\frac {5}{21}}&0&1/3&-
1/3\\\noalign{\medskip}1/7&1/7&1/7&2/7&2/7&0&-1&1\\\noalign{\medskip}0
&0&0&0&0&0&0&1\\\noalign{\medskip}-1/7&-1/7&-1/7&-2/7&-2/7&1&1&-1
\end {array}\right ]$$

We have therefore the Mellin transform of the period integral 
associated to the CI $\{x $ $\in$ $(\bC
^\times)^5$ $;$ $f_1(x)$ $+$ $s_1 = 0, f_2(x) + 1 = 0\}$ 
up to 
$7-$periodic functions,
$$ M_{0,\gamma}^{0}(z)=
\Gamma( -\frac{1}{7}(z_1-1))^3
\Gamma(-\frac{2}{7}(z_1-1))^2 \Gamma(z_1) =
\frac{\Gamma(\xi^{(1)})^3\Gamma(2\xi^{(1)})^2}{\Gamma(7\xi^{(1)})}.$$
After the construction $(1.4)^T$ we see that
$$ \;^Tf_1(x)= x_1^7 + x_2^7+x_3^7+ x_2 x_4^3+x_3x_5^3,\;\;
\;^Tf_2=x_1x_2x_3x_4x_5.$$

We then have the Mellin transform of the period integral associated to the
CI $\{x \in (\bC^\times)^5; \;^Tf_1 ( x ) + s_1 = 0 , \;^Tf_2(x) + 1 = 0\}$
up to $147-$periodic functions,
$$ M_{0,\gamma}^{0}(z)=
\Gamma( -\frac{1}{7}(z_1-1))
\Gamma(-\frac{2}{21}(z_1-1))^2
\Gamma(-\frac{1}{3}(z_1-1))^2 \Gamma(z_1) =
\frac{\Gamma(3\xi^{(1)})\Gamma(2\xi^{(1)})^2
\Gamma(7\xi^{(1)})^2}{\Gamma(21\xi^{(1)})}.$$


\vspace{\fill}

%

\noindent

\begin{flushleft}
\begin{minipage}[t]{6.2cm}
  \begin{center}
{\footnotesize Indepent University of Moscow\\
Bol'shoj Vlasievskij pereulok 11,\\
 Moscow, 121002,\\
Russia\\
{\it E-mails}:  tanabe@mccme.ru, tanabe@mpim-bonn.mpg.de \\}
\end{center}
\end{minipage}
\end{flushleft}

\end{document}